\documentclass[reqno]{amsart} 
\usepackage{hyperref}
\usepackage{amsmath,amsthm,amsfonts,amscd,amssymb,amsbsy,epsf}
\usepackage{verbatim}
\usepackage{graphicx}
\usepackage{enumerate}
\usepackage{bbm, dsfont, mathrsfs}
\DeclareMathOperator*{\essinf}{ess\,inf}
\usepackage{xcolor}

\def\V{\mathbf V}
\def\g{\mathbf g}
\def\n{\boldsymbol{ \widehat {n} } }
\def\outer{\boldsymbol{\widehat{\nu} } }
\def\outerth{\boldsymbol{\widetilde{\nu} } }
\def\u{\mathbf u}
\def\v{\mathbf v}
\def\w{\mathbf w}
\def\x{\mathbf x}
\def\y{\mathbf y}
\def\div{\boldsymbol{\nabla } \cdot }
\def\grad{\boldsymbol{\nabla}}
\def\gradth{\boldsymbol{\widetilde{\nabla} } }

\def\ind{\boldsymbol{\mathbbm {1} } }
\def\eversor{\boldsymbol{ \widehat {e}} }
\def\defining{\overset{\mathbf{def}}=}

\def\uone{\mathbf{u}^{1}}
\def\utwo{\mathbf{u}^{2}}
\def\ueps{\mathbf{u}^{\epsilon}}
\def\unormeps{u^{\,\epsilon,\,2}_{\,3}}

\def\utaneps{\mathbf{\widetilde{u}}^{\,\epsilon,\,2}}

\def\utangeps{\mathbf{u}^{\,\epsilon,\,2}_{\,\tau}}
\def\peps{p^{ \epsilon}}
\def\unorm{u_{3}^{\,2}}
\def\utan{\mathbf{\widetilde{u}}^{\,2}}
\def\utang{\mathbf{u}_{\,\tau}^{\,2}}
\def\pepsone{p^{ \epsilon, 1}}
\def\pepstwo{p^{ \epsilon, 2} }
\def\uepsone{\mathbf{u}^{ \epsilon, 1} }
\def\uepstwo{\mathbf{u}^{ \epsilon, 2} }

\def\R{\mathbbm{R} } 
\def\Hdiv{\mathbf{H_{ div } } }
\def\Ltwo{\mathbf{L^{\! 2} } }
\def\wone{\mathbf{w}^{ 1}}

\def\wtang{\mathbf{w}_{\,\tau}^{\,2}}

\def\vone{\mathbf{v}^{1}}
\def\vtwo{\mathbf{v}^{2}}
\def\vtan{\mathbf{\widetilde{v}}}
\def\vtang{\mathbf{v}_{\tau}^{2}}
\def\vnorm{v_{3}}
\def\qone{q^{ 1}}
\def\qtwo{q^{ 2}}

\def\pone{p^{1}}
\def\ptwo{p^{2}}
\def\aone{a_{1}}
\def\atwo{a_{2}}

\def\MTtang{M^{\,T, \boldsymbol{\tau} } }
\def\MTnormal{M^{\,T, \n } }
\def\MNtang{M^{\kversor, \boldsymbol{\tau}}}
\def\MNnormal{M^{\kversor, \n } }
\def\iversor{\boldsymbol {\widehat{i } } }
\def\jversor{\boldsymbol {\widehat{j} } }
\def\kversor{\boldsymbol{\widehat {k} } }
\def\eone{\boldsymbol {\widehat {e} }_1}
\def\etwo{\boldsymbol {\widehat {e} }_2}
\def\xthilde{\widetilde{\mathbf{x}}}
\def\xthree{x_{3}}

\def\wtang{\mathbf{w}_{\,\tau}^{\,2}}

\def\iindex{1\,\leq\, i\, \leq\, I}
\def\jindex{0\,\leq\, j\, \leq\, I}

\def\Gi{G_{i}}
\def\zetai{\zeta_{\,i}}
\def\hi{h_{i}}
\def\hj{h_{\,j}}
\def\Omegai{\Omega_{i}}
\def\Omegaj{\Omega_{j}}
\def\Gammai{\Gamma_{i}}

\def\Gammab{\Gamma_{\mathbf{b}}}
\def\Gammat{\Gamma_{\mathbf{t}}}
\def\ni{\boldsymbol{\widehat{n}}^{ (i) } }

\def\Thetaj{\Theta_{j}}
\def\Thetai{\Theta_{i}}
\def\lambdai{\lambda_{\,i}}
\def\Lambdai{\Lambda_{\,i}}
\def\Lambdaip{\Lambda_{\,i}^{+}}
\def\Lambdain{\Lambda_{\,i}^{-}}

\begin{document}

\title[Geological Fissured Systems]
{On the Homogenization of Geological Fissured Systems With Curved non-periodic Cracks}

\author[Fernando A. Morales]
{Fernando A. Morales}

\address {Fernando A. Morales\newline 
Escuela de Matem\'aticas, Universidad Nacional de Colombia, Sede Medell\'in. Colombia}
\email{famoralesj@unal.edu.co}


\thanks{The author was supported by projects HERMES 14917, HERMES 17194 from Universidad Nacional de Colombia, Sede Medell\'in, and by the Department of Energy, Office of Science, USA through grant 98089.}

\subjclass[2000]{35F15, 80M40, 76S99, 35B25}

\keywords {fissured media, tangential flow, interface geometry, coupled Darcy flow system, upscaling, mixed formulations}

\begin{abstract}
We analyze the steady fluid flow in a porous medium containing a network of thin fissures i.e. width $\mathcal{O}(\epsilon)$, where all the cracks are generated by the rigid translation of a continuous piecewise $C^{1}$ functions in a fixed direction. The phenomenon is modeled in mixed variational formulation, using the stationary Darcy's law and setting coefficients of low resistance $\mathcal{O}(\epsilon)$ on the network. The singularities are removed performing asymptotic analysis as $\epsilon \rightarrow 0$ which yields an analogous system hosting only tangential flow in the fissures. Finally the fissures are collapsed into two dimensional manifolds. 
\end{abstract}
\maketitle
\allowdisplaybreaks
\newtheorem{theorem}{Theorem}[section]
\numberwithin{equation}{section}
\newtheorem{lemma}{Lemma}[section]
\newtheorem{corollary}[lemma]{Corollary}
\newtheorem{definition}[theorem]{Definition}
\newtheorem{proposition}[theorem]{Proposition}
\newtheorem{remark}[lemma]{Remark}
\newtheorem{hypothesis}[theorem]{Hypothesis}

\section{Introduction}   
Groundwater and oil reservoirs are frequently fissured or layered i.e. the bed rock contains fissures of characteristic dimensions considerably higher than those of the average pore size of the rock. The modeling of saturated flow through geological structures such as these, gives rise to singular problems of partial differential equations \cite{Show97}. On one hand the singularities are due to the drastic change of permeability from the rock matrix to the fissures. On the other hand a geometric singularity is introduced due to the thinness of the fractures. The presence of singularities in the model has non-desirable effects in their numerical implementation; some of these are ill-condition matrices, high computational costs, numerical stability, etc. This subject is a very active research field, see \cite{ArbBrunson2007, Gatica2009, Yotov, Loredana1, JaffRob05} for numerical analysis aspects, \cite{ZhaoQin, Levy83} for modeling discussion and \cite{Allaire2009, ArbLehr2006, Gunzburger2009, Mikelic89, MoralesShow1} for rigorous mathematical treatment of the phenomenon. Homogenization and asymptotic analysis techniques are a common approach for the analytical point of view. However, the remarkable achievements in the field require very restrictive hypotheses for the description of the geometry such as uniformly distributed, regular geometric shapes or periodic arrayed structures \cite{Hornung, SP80}. In general the variational methods for partial differential equations can formulate successfully a wide class of geometric domains, the limited treatment of the geometry comes from the notorious difficulties it introduces in the asymptotic analysis of the problem. 

In the present work, the geometric possibilities of the medium are broaden to an unprecedented setting: free from the aforementioned hypotheses. We use the \emph{mixed mixed formulation} and the scaling for the flow resistance coefficients presented in \cite{MoralesShow2}, then a careful choice of directions or ``stream lines'', consistent with the natural scaling of the problem permits a successful asymptotic analysis of the model. This leads to a system coupled though multiple two dimensional manifolds representing the fissures in the upscaled model. Additionally, the formulation allows remarkable generality in the fluid exchange balance conditions between the rock matrix and the channels, substantial efficiency for handling the system of equations as well as the information (coefficients, matrices, etc) describing the geometry of the fractures, mostly due to the fact that it does not demand coupling constraints on the underlying spaces of functions. The main goal of the paper is to emphasize on the geometry, consequently the study is limited to the steady case. We describe flow with Darcy's law 
\begin{subequations}\label{Pblm FORMAL fissures system strong problem}
\begin{equation}\label{Eq FORMAL Darcy}
a (\cdot)\,\u + \grad p + \g = 0 ,
\end{equation}
together with the conservation law
\begin{equation}\label{Eq FORMAL fluid mass conservation}
\div \u = F.
\end{equation}
Drained and non-flux boundary conditions on different parts of the domain boundary will be specified to set a boundary value problem. The fluid exchange across the interface separating the regions are given by
\begin{equation}\label{Eq FORMAL interface balance 1 gamma}
\pone - \ptwo = \alpha \, \uone\quad \text{and}
\end{equation}
\begin{equation}\label{Eq FORMAL interface balance 2 gamma}
\uone \cdot \n - \utwo\cdot\n  = 
f_{\scriptscriptstyle \Gamma} 
\quad \text{on}\;\, \Gamma .
\end{equation}
\end{subequations}
Here, the coefficient $a(\cdot)$ is the flow resistance i.e. the fluid viscosity times the inverse of the permeability of the medium, to be scaled consistently with the fast and slow flow regions of the medium. Finally, the coefficient $\alpha$ indicates the fluid entry resistance of the rock matrix. 

In the following section we define the geometric setting, formulate the problem in mixed mixed variational formulation and establish its well-posedness. In section three the problem is referred to a common geometric setting in order make possible the asymptotic analysis, the existence of a-priori estimates and the structure of the limiting solution are also shown. Section four studies the formulation and well-posedness of the limiting problem and finds its strong form, particularly important for boundary and interface conditions and proves the strong convergence of the solutions. Section five sets the limiting problem as a coupled system with two dimensional interfaces and section six discusses the possibilities and limitations of the technique as well as related future work.        
%
%
%
%
%
%
%
%
%
%
\section{Formulation and Geometric Setting}   
Vectors are denoted by boldface letters as are vector-valued functions and corresponding function spaces. We use $\xthilde$ to indicate a vector in $\R^{\! 2}$; if $\x\in \R^{\! 3}$ then the $\R^{\! 2}\times\{0\}$ projection is identified with $\xthilde\defining(x_{1}, x_{2})$ so that $\x = (\xthilde, \xthree)$. The symbol $\gradth$ represents the gradient in the first first two directions: $\iversor$, $\jversor$. Given a function $f: \R^{\! 3}\rightarrow \R$ then  $\int_{\mathscr {M} } f\,dS$ is the notation for its surface integral on the $\R^{\! 2}$ manifold $\mathscr{M}\subseteq \R^{\! 3}$. $\int_{A} f\, d\x$ stands for the volume integral in the set $A\subseteq \R^{\! 3}$; whenever the context is clear we simply write $\int_{A} f$. In the same fashion, whenever there is no confusion $\sum_{\,i}$, $\prod_{\,i}$ indicate $\sum_{i = 1}^{I}$ and $\prod_{i = 1}^{I}$ respectively. 
\begin{equation}\label{Def general vertical shift}
A+t\defining\left\{\x+ t\,\kversor: \x\in A\right\}
\end{equation}
The symbol $\boldsymbol{\widehat{\nu} }$ denotes the outwards normal vector on the boundary  of a given domain $\mathcal{O}$ and $\n$ denotes the normal upwards vector to a given surface i.e. $\n \cdot \kversor \geq 0$. For any $A\subseteq\R^{3}$ and $t\in \R$ we define its $t$-vertical shift by
%
%
\subsection{General Geometric Setting}   
The present work will be limited to the study of fractured media where each fissures can be described in a specific way. 
\begin{definition}\label{Def eligible surface and generated fissure}
Let $G\subseteq \R^{2}$ be open a bounded open simply connected set and $\zeta\in C(\overline{G})$ be a piecewise $C^{1}$. Define the surface
\begin{equation}\label{Def C2 surface}
\Gamma\defining \left\{\left[\xthilde, \zeta\left(\xthilde\right)\right]:\xthilde\in G\right\} .
\end{equation}
We say $\Gamma$ is a surface eligible for vertical translation fissure generation if $\essinf\{ \n(s)\cdot\kversor: s\in \Gamma\}>0$. Given vertical height $h>0$ define the fissure of height $h$ generated by a rigid vertical translation of $\Gamma$ by the domain
\begin{equation}\label{Def vertical fissure}
\Omega\left(h, \Gamma\right)\defining \left\{\left(\xthilde, y\right):
\zeta\left(\xthilde\right)< y < \zeta\left(\xthilde\right)+h\right\} .
\end{equation}
\end{definition}
\begin{remark}\label{Rem comments on the height width relation}
Notice that in the definition of $\Omega(h, \Gamma)$ we mention $h$ as the height and not as the width of the crack. Figure \eqref{Fig Fissure with high gradients} shows that, depending on the gradient of the surface the height $h$ can become significantly different from the actual width. 
\end{remark}
\begin{figure}[!]
\caption[1]{Unidirectional Translation Generated Fissures}\label{Fig Fractured Medimu}
\centerline{\resizebox{10cm}{7cm}
{\includegraphics{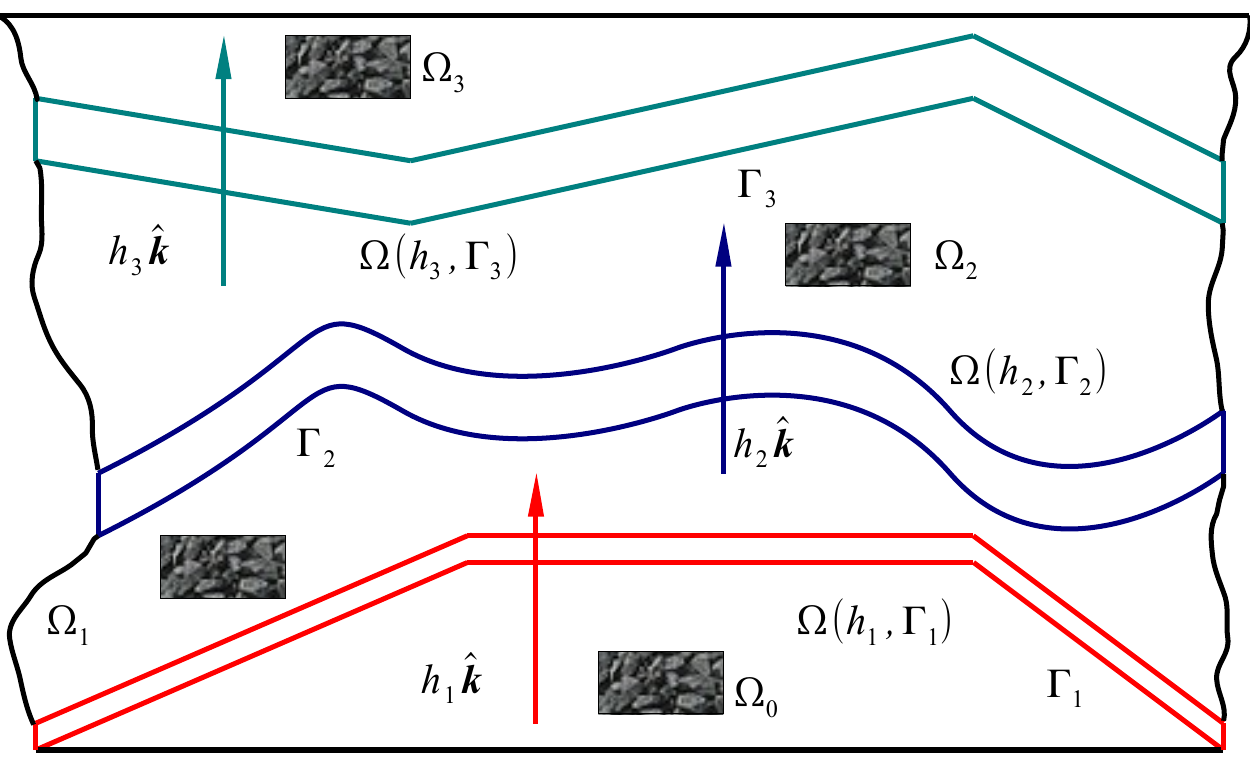} } }
\end{figure}
The analysis will be limited to the type of geological system shown in figure \eqref{Fig Fractured Medimu}. It depicts a region $\Omega\subseteq\R^{\!3}$ containing a network of fissures generated by vertical rigid translation continuous piecewise $C^{1}$ surfaces. Such a region is completely characterized in the following definition
\begin{definition}\label{Def unidirectional translation medium}
We say a totally fractured medium of vertical translation generated fissures is a finite collection of
\begin{subequations}\label{Def information of the domain}

Surface functions
\begin{multline}\label{Def fissures surfaces}
\{\zetai\in C (\overline{\Gi}):\Gi\subseteq\R^{2}\;\text{open bounded simply connected region} ;\\
\zetai \,\;\text{piecewise}\;\,C^{1} \,\;\text{functions such that}\,\; \essinf \ni\cdot\kversor >0 
,\, 1 \leq i \leq I\}.
\end{multline}

vertical heights
\begin{equation}\label{Def scalars width}
\left\{\hi >0 :\iindex\right\} ,
\end{equation}

and rock-matrix regions
\begin{equation}\label{Def rock matrix}
\left\{\Omegai \subseteq \R^{3}: \Omegai \neq \emptyset\;\text{open bounded simply connected region},\; 0 \leq i\leq I\right\} .
\end{equation}
%
%
\end{subequations}
Verifying the following properties:
\begin{subequations}\label{Def properties of the medium}

Non-overlapping condition and indexed ordered
\begin{equation}\label{Def ordering condition 1}
\sup \left\{\zetai(\xthilde) + h_{i}: \xthilde\in \Gi\right\} <
\inf \left\{\zeta_{\,i+1}(\xthilde): \xthilde\in G_{i+1}\right\}
\;\forall\; 1\leq i\leq I-1 .
\end{equation}

The interface-domain condition
\begin{equation}\label{Def interface boundary conditions}
\begin{split}
\partial \Omegai \cap \partial\Omega(h_{i+1}, \Gamma_{i+1}) = \Gamma_{i+1}
\quad \forall\; 0\leq i\leq I-1 ,\\
\partial \Omegai \cap \partial\Omega(\hi, \Gammai)= \Gammai+\hi
\quad \forall \,\iindex .
\end{split}
\end{equation}

And the condition of connectivity only through fissures 
\begin{equation}\label{Def fissures connectedness condition}
cl(\Omega_{\,\ell})\cap cl(\Omega_{\,k}) = \emptyset\,,\quad\text{whenever}\;\ell\neq k .
\end{equation}
\end{subequations}

For convenience of notation define $\Gamma_{0} \defining \partial \Omega_{0} - \Gamma_{1}$ and $h_{0}\defining 0$. The fissured system described above will be denoted $\left\{\left(\Gammai, \hi, \Omegai\right):0 \leq i\leq I\right\}$. The sets $\Omega_{1}, \Omega_{2}$ are the rock matrix and the fissures regions respectively \emph{i.e.}
\begin{equation}\label{Def fissured medium region}
\begin{split}
\Omega_{1} \defining \bigcup_{i = 0}^{I}\Omegai\,&,\quad
\Omega_{2}\defining  \bigcup_{i = 1}^{I} \Omega\left(\hi, \Gammai\right)\\
\Omega & \defining  \Omega_{1}\cup\Omega_{2}.
\end{split}
\end{equation}
The global bottom and top interfaces are defined by
\begin{equation}\label{Def global interfaces}
\begin{split}
\Gammat \defining \bigcup_{i = 1}^{I}\Gammai\,&,\quad
\Gammab \defining \bigcup_{i = 1}^{I} \Gammai+\hi\\
\Gamma &\defining \Gammab\cup\Gammat
\end{split}
\end{equation}
Finally, $\ni$ indicates the upwards normal vector to the surface $\Gammai$ \emph{i.e.}
\begin{equation}\label{Def Standarazied upwards normal vector}
   \ni \defining \frac{(-\gradth \zetai, 1)}{\vert (-\gradth \zetai, 1)\vert} .
\end{equation}
When there is no confusion $\n$ denotes the normal vector with respect to the surface of the crack.
\end{definition}
\begin{remark}\label{Rem connectivity condition on the domains}
The condition \eqref{Def fissures connectedness condition} of connectivity only through fissures is not required for modeling the problem in mixed formulation as it is presented in section \eqref{Sec mixed formulation}; however it is necessary for the asymptotic analysis of the system. The same holds for the requirement of simply connected domains. 
\end{remark}
%
%
%
%
\subsection{A Local System of Coordinates}\label{Sec local system of coordinates}
Some aspects of the flow through the fissures are handled more conveniently when the velocities are expressed in a coordinate system consistent with the geometry of the surface that generates the crack. Let $\Gamma$ be a surface as defined in \eqref{Def eligible surface and generated fissure} and $\n$ the upwards normal to the surface $\Gamma$ i.e. $\n = \n(s) = \n (\xthilde)$. Now, for each point $\xthilde$ we choose a local orthonormal basis in the following way
\begin{equation}\label{Eq local orthogonal basis}
\mathcal{B}(\xthilde)\defining
\left\{\eone(\xthilde), \etwo (\xthilde), \n(\xthilde)\right\}
\end{equation}
Let $M = M (\xthilde)$ be the orthogonal matrix relating the global canonical basis with the local one i.e.
\begin{subequations}\label{Def orthogonal matrix}
\begin{equation}\label{Def orthogonal matrix 1 }
M (\xthilde) \, \iversor = \eone(\xthilde)
\end{equation}
\begin{equation}\label{Def orthogonal matrix 2}
M (\xthilde) \, \jversor = \etwo (\xthilde)
\end{equation}
\begin{equation}\label{Def orthogonal matrix 3}
M (\xthilde) \, \kversor = \n (\xthilde) 
\end{equation}
\end{subequations}
The block matrix notation for this local matrix will be
\begin{equation}\label{Eq blocks of the rotational matrix}
M (\xthilde) \defining
\left(\begin{array}{cc}
\MTtang & \MTnormal\\[7pt]
\MNtang & \MNnormal
\end{array}\right)(\xthilde)
\end{equation}
Here the index $T$ stands for the first two components in the directions $\iversor, \, \jversor$ while the index $\tau$ stands for the expression of the velocity orthogonal to the component in the direction $\n$. Then $\w =  [\w_{\tau}, \w_{\n}](\xthilde)$ with the following relations
\begin{subequations}\label{Eq definition of normal and tangential velocity}
\begin{equation}\label{Eq definition of normal velocity}
w_{\n} \defining \w\cdot\n (\xthilde)
\end{equation}
\begin{equation}\label{Eq definition of tangential velocity}
\w_{\tau} \defining \left(\,\w\cdot \eone (\xthilde), \w\cdot\etwo(\xthilde)\,\right)
\end{equation}
\end{subequations}
Clearly, the relationship between velocities is given by
\begin{multline}\label{Eq local and global velocities}
\w \left(\xthilde, x_{\,3}\right)
= \left\{\begin{array}{c}
\widetilde{\w}\\[7pt]
\w\cdot\kversor 
\end{array}
\right\}
\left(\xthilde, x_{\,3}\right)
= M (\xthilde) \left\{\begin{array}{c}
\w_{\tau}\\[7pt]
\w\cdot\n
\end{array}
\right\}
\left(\xthilde, x_{\,3}\right)\\[5pt]
= \left(\begin{array}{cc}
\MTtang (\xthilde) & \MTnormal (\xthilde)\\[7pt]
\MNtang (\xthilde) & \MNnormal (\xthilde)
\end{array}\right)
\left\{\begin{array}{c}
\w_{\tau}\\[7pt]
\w\cdot\n
\end{array}
\right\}
\left(\xthilde, x_{\,3}\right)
\end{multline}
\begin{proposition}\label{Th isometry change of coordinates}
Let $h>0$, $\Gamma$, $\Omega(h, \Gamma)$ be as in definition \eqref{Def eligible surface and generated fissure}; $\n$ be the upwards normal to the surface $\Gamma$ and $M$ be the matrix defined by \eqref{Def orthogonal matrix}. Then
\begin{enumerate}[(i)]
   \item 
   The map $\w\mapsto M(\xthilde) \w$ is an isometry in $\Ltwo(\Omega(h, \Gamma) )$. In particular if $\w_{\tau}, \w\cdot\n$ are defined as in \eqref{Eq definition of normal and tangential velocity} then $\w\in \Ltwo(\Omega(h, \Gamma) )$ if and only if $\w_{\tau}\in L^{2}(\Omega(h, \Gamma) )\times L^{2}(\Omega(h, \Gamma))$ and $\w\cdot \n \in L^{2}(\Omega(h, \Gamma))$.
   
   \item
   If $\w\in \Ltwo(\Omega(h, \Gamma))$ is such that $\partial_{z} \w \in \Ltwo(\Omega(h, \Gamma))$ then
\begin{equation}\label{Eq partial z velocities}
\partial_z \w \left(\xthilde, z\right)
= M (\xthilde) \left\{\begin{array}{c}
\partial_z\,\w_{\tau}\\[7pt]
(\partial_z\,\w)\cdot\n
\end{array}
\right\}
\left(\xthilde, z\right)
\end{equation}
   \end{enumerate}
\begin{proof}
\begin{enumerate}[(i)]
   \item 
   For $\xthilde$ fixed the matrix $M(\xthilde)$ is orthogonal i.e. for arbitrary functions $\v,\w \in \Ltwo(\Omega(h, \Gamma) )$ and $\x\in \Omega(h, \Gamma)$ holds $\v(\x)\cdot\w(\x) = M(\xthilde) \v(\x)\cdot M(\xthilde) \w(\x)$. Hence  
   \begin{multline*}
   \int_{\Omega(h, \Gamma) }\v(\x)\cdot\w(\x)\,d\x =
   \int_{\Omega(h, \Gamma) } M(\xthilde)\v(\x)\cdot M(\xthilde)\w(\x)\,d\x \\
   =\int_{\Omega(h, \Gamma) } \v_{\tau}(\x)\cdot \w_{\tau}(\x)\,d\x 
   + \int_{\Omega(h, \Gamma) } (\v\cdot\n) (\x)\cdot (\w\cdot\n)(\x)\,d\x
   \end{multline*}
   The equality of the second line shows the necessity and sufficiency of the tangential and normal components been square integrable in the domain $\Omega(h, \Gamma)$.
   
   \item
   It follows from a direct calculation of distributions with $\boldsymbol {\varphi}\in [C_{0}^{\infty} (\Omega(h, \Gamma) ) ]^{3}$ arbitrary and the fact that $\partial_{z} M = 0$.
   
   \end{enumerate}
\end{proof}
\end{proposition}
%
%
%
%
\subsection{The Problem and its Formulation}\label{Sec strong general problem}
In this section we define the problem in a rigorous way and give a variational formulation in which it is well-posed. Let $\left\{\left(\Gammai, \hi, \Omegai\right):\iindex\right\}$ be a totally fractured domain of vertical translation generated fissures. We denote $\vone, \pone$ the velocity and pressure in the rock matrix region  $\Omega_{1}$. In the same fashion $\vtwo, \ptwo$ denote the velocity and pressure in the fissures region $\Omega_{\,2}$.
Consider the problem
\begin{subequations}\label{Pblm fissures system strong problem}
\begin{equation}\label{Eq Darcy rock matrix}
\aone \uone + \grad \pone + \g = 0\quad \text{and}
\end{equation}
\begin{equation}\label{Eq fluid mass conservation rock matrix}
\div \uone = F\quad \text{in}\;\Omega_{1} .
\end{equation}
\begin{equation}\label{Eq drained condition rock matrix}
\pone = 0 \quad \text{on} \; \partial\Omega_{1}-\Gamma .
\end{equation}
\begin{equation}\label{Eq interface balance 1 gamma}
\pone - \ptwo = \alpha \, \uone\cdot\n \; \ind_{ \Gammab } 
-\alpha \, \uone\cdot\n \; \ind_{ \Gammat }
\quad \text{and}
\end{equation}
\begin{equation}\label{Eq interface balance 2 gamma}
\left(\uone -\utwo \right)\cdot\n\;\ind_{ \Gammat } 
-\left(\uone -\utwo \right)\cdot\n\; \ind_{ \Gammab } = 
f_{\scriptscriptstyle \Gamma} 
\quad \text{on}\; \Gamma .
\end{equation}
\begin{equation}\label{Eq Darcy fissures region}
\atwo\,\utwo+\grad\ptwo +\g = 0\quad \text{and}
\end{equation}
\begin{equation}\label{Eq fluid mass conservation fissures}
\div \utwo = F \quad \text{in} \;\Omega_{2} .
\end{equation}
\begin{equation}\label{Eq non flux condition fissures}
\utwo \cdot\n = 0 \quad \text{on}\, \partial \Omega_{2}-\Gamma .
\end{equation}
\end{subequations}
The flow resistance coefficients $\aone, \atwo$ and the fluid entry resistance coefficient $\alpha$ are assumed to be positively bounded from below and above, see \cite{MoralesShow2}. In equations \eqref{Eq interface balance 1 gamma}, \eqref{Eq interface balance 2 gamma} the split of cases is made in order to be consistent with the sign of the upwards normal vector $\n$. 
%
%
%
%
%
%
%
%
%
%
\subsection{Mixed Formulation of the Problem}\label{Sec mixed formulation}
We start defining the spaces of velocities and pressures
\begin{subequations}\label{Def spaces of functions}
\begin{equation}\label{Def space of velocities}
\V\defining\{\v\in\Ltwo(\Omega): \div\vone\in \Ltwo(\Omega_{1})
\,,\vone\cdot\n\,\vert_{\Gamma}\in L^{2}(\Gamma)\} .
\end{equation}
\begin{equation}\label{Def space of pressures}
Q\defining\{q\in L^{2}(\Omega): \grad\qtwo\in \Ltwo(\Omega_{\,2}) \}
\end{equation}
Endowed with their natural norms
\begin{equation}\label{Def normm space of velocities}
\Vert\, \v\,\Vert_{\scriptscriptstyle \V} \defining \{\,\Vert \,\v\,\Vert_{\Ltwo(\Omega)}^{2}+\Vert\, \div\vone\,\Vert_{L^{2}(\Omega_1)}^{2} + \Vert\, \vone\cdot\n\,\Vert_{L^{2}(\Gamma)}^{2}\}^{1/2}
\end{equation}
\begin{equation}\label{Def norm of pressures}
\Vert\, q\,\Vert_{\scriptscriptstyle Q} \defining \{\,\Vert\, q\,\Vert_{L^{2}(\Omega)}^{2}
+\Vert \,\grad\qtwo\,\Vert_{L^{2}(\Omega_{\,2})}^{2}\,\}^{1/2}
\end{equation}
\end{subequations}
\begin{remark}\label{Rem notation simplification of Gamma}
   In the spaces above it is understood that 
   \begin{multline}\label{Eq notation condensation}
   \Vert \v\cdot \n\Vert_{\scriptscriptstyle L^{2}(\Gamma)}^{2} = 
   \Vert \v\cdot \n\Vert_{\scriptscriptstyle L^{2}(\Gammab)}^{2} +
   \Vert \v\cdot \n\Vert_{\scriptscriptstyle L^{2}(\Gammat)}^{2} \\
   = \sum_{\,i = 1}^{I} \Vert \v\cdot \ni\Vert_{\scriptscriptstyle L^{2}(\Gammai)}^{2} + 
   \sum_{\,i = 1}^{I} \Vert \v\cdot \ni\Vert_{\scriptscriptstyle L^{2}(\Gammai + \hi)}^{2}
   \end{multline}
\end{remark}
Consider the problem
\begin{subequations}\label{Pblm weak form epsilon free system}
\begin{equation*}
\text {Find } \; p\in Q, \, \u\in \V
\end{equation*}
\begin{multline}\label{Pblm weak form epsilon free system 1}
\int_{\Omega_1} a_1 \,\u \cdot \v 
+  \int_{\Omega_2} a_2\, \u \cdot \v
- \int_{\Omega_1} p \,\div\v  
+ \int_{\Omega_2} \grad p \cdot\v 
\\
+\alpha\int_{\Gamma}\left(\uone\cdot\n\right)\left(\vone\cdot\n\right) d S
%
-\int_{\,\Gammat}\ptwo \left(\vone\cdot\n\right)\, d S
+\int_{\,\Gammab}\ptwo \left(\vone\cdot\n\right)\, d S
= - \int_{\Omega}\g \cdot \v
\end{multline}
\begin{multline}\label{Pblm weak form epsilon free system 2}
\int_{\Omega_1}\div\u\, q
- \int_{\Omega_{2}} \u \cdot \grad q\\
+ \int_{\,\Gammat}\left(\uone\cdot\n\right) \qtwo\, d S
- \int_{\,\Gammab}\left(\uone\cdot\n\right) \qtwo\, d S
= \int_{\Omega}F\, q
+ \int_{\Gamma} f_{\scriptscriptstyle \Gamma}\, \qtwo\,d S
\end{multline}
\begin{equation*}
\text{ for all }  q\in Q,\,\v\in\V
\end{equation*}
\end{subequations}
\begin{remark}\label{Rem sign explanation}
In the formulation above the non-symmetric interface terms are split in two pieces in order to express everything in terms of the upwards normal vector $\n$. In the case of the symmetric term $\int_{\Gamma}(\uone\cdot\n)(\vone\cdot\n)dS$ in \eqref{Pblm weak form epsilon free system 1} such split becomes unnecessary since the sign of the normal vector changes in both factors canceling each other.
\end{remark}
Define the bilinear forms $\mathcal{A}:\V\rightarrow\V'$, $\mathcal{B}:V\rightarrow Q'$, $\mathcal{C} :Q\rightarrow Q'$ by
\begin{subequations}\label{Def operators fissures system}
\begin{equation}\label{Def velocity operator fissures system}
\mathcal{A}\v(\w)\defining \int_{\Omega_{1}}\aone\,\v\cdot\w
+\int_{\Omega_{2}}\atwo\,\v\cdot\w +
\alpha \int_{\Gamma}\left(\vone\cdot\n\right)\left(\wone\cdot\n\right)\,d S
\end{equation}
\begin{equation}\label{Def mixed operator fissures system}
\mathcal{B}\v (q) \defining -\int_{\Omega_{1}}\div\v\,q
+\int_{\Omega_{2}}\v\cdot \grad q
-\int_{\Gammat} (\vone\cdot\n) \,\qtwo \, d S
+\int_{\Gammab} (\vone\cdot\n) \,\qtwo \, d S
\end{equation}
\end{subequations}
Then, the system \eqref{Pblm weak form epsilon free system} is a mixed formulation for the problem \eqref{Pblm fissures system strong problem} with the abstract form
\begin{equation}\label{Pblm mixed formulation}
\begin{split}
\u\in \V, \;p\in Q  : \quad \mathcal{A}\u+\mathcal{B}\,'p &= -\mathbf{g}\;\text{in}\;\V ' ,\\
-\mathcal{B}\,\u &= f\quad\text{in}\;Q ' .
\end{split}
\end{equation}
For the sake of completeness recall some well known results
\begin{theorem}\label{Th General Mixed Formulation well posedness}
   Let $\V, Q$ be Hilbert spaces and $\Vert \cdot \Vert_{\V}, \Vert \cdot \Vert_{Q}$ be their respective norms. Let $\mathcal {A}: \V \rightarrow \V '$, $\mathcal {B}: \V \rightarrow Q'$ be continuous linear operators such that
   \begin{enumerate}[(i)]
      \item 
      $\mathcal {A}$ is non-negative and $\V$-coercive on $\ker \mathcal{B}$. 
      
      \item
      The operator $\mathcal{B}$ satisfies the inf-sup condition
      \begin{equation}\label{Ineq general inf sup condition}      
         \inf_{q\in Q} \, \sup_{\v\in \V} \frac{\vert \mathcal{ B } \v (q) \vert }{\Vert \v\Vert_{\scriptscriptstyle \V} \, \Vert q\Vert_{\scriptscriptstyle Q}}  > 0 .
      \end{equation}
      Then, for each $\g\in \V'$ and $f\in Q'$ there exists a unique solution $[\u, p]\in \V\times Q$ to the problem \eqref{Pblm mixed formulation}. Moreover, it satisfies the estimate
      \begin{equation}\label{Ineq well posedness estimate}     
         \Vert \u\Vert_{\scriptscriptstyle \V} + \Vert p \Vert_{\scriptscriptstyle Q}
         \leq K \left (\Vert \g\Vert_{\scriptscriptstyle \V'} + \Vert f\Vert_{\scriptscriptstyle Q'}\right)
      \end{equation}
   \end{enumerate}
   \begin{proof}
      See \cite{GiraultRaviartFEM}
   \end{proof}
\end{theorem}
\begin{lemma}\label{Th control with trace and gradient}
Let $\mathcal{O}$ be an open connected bounded set in $\R^{\! N}$ and $\mathcal{G}\subseteq \partial\mathcal{O}$ with non-null $\R^{\! N-1}$-Lebesgue measure, then there exists $\kappa = \kappa(\mathcal{O}) >0$ such that
\begin{equation} \label{Eq control with trace and gradient}
\Vert \grad \eta \Vert_{\Ltwo(\mathcal{O})} + \Vert  \eta \Vert_{\mathcal{G}}\geq \kappa\,
\Vert  \eta \Vert_{H^{1}(\mathcal{O})}
\end{equation}
for all $\eta\in H^{1}(\mathcal{O})$.
   \begin{proof}
      See proposition 5.2 of \cite{Showalter97} or lemma 1.2 in \cite{MoralesShow2}.
   \end{proof}
\end{lemma}
\begin{corollary}\label{Th control by gradient and trace}
   There exists a constant $\kappa >0$ such that 
   \begin{equation}\label{Ineq control by gradient and trace}
      \Vert \grad q \Vert_{\Ltwo(\Omega_{2})} ^{2} 
      + \Vert  q \Vert_{L^{2}(\Gamma)} ^{2}
      \geq\kappa \,
      \Vert  q \Vert_{L^{2}(\Omega_{2})} ^{2} .
   \end{equation}
   For all $q\in H^{1}(\Omega_{2})$
   \begin{proof}
      Apply lemma \eqref{Th control with trace and gradient} on each connected component $\Omega(\hi, \Gammai)$ and choose $\kappa$ as the minimum constant associated to each domain. 
   \end{proof}
\end{corollary}
\begin{lemma}\label{Th inf-sup condition original}
   The operator $\mathcal {B}$ satisfies the inf-sup condition \eqref{Ineq general inf sup condition}. 
   \begin{proof}
      We use the same strategy presented lemma 1.3 in \cite{MoralesShow2} with a slight modification in the construction of the particular test function. Fix $q\in Q$ and denote $\xi_{j}$ the unique solution of the problem
      \begin{equation}\label{Pblm Test Function}
      \begin{split}
         -\div \grad \xi_{j} = \qone \; \text{in} \; \Omega_{j}, \\
         \grad \xi_{j} \cdot \n = \qtwo \; \text{on} \; \Gamma_{j} ,
         \quad \grad \xi_{j} \cdot \n = - \qtwo \; \text{on} \; \Gamma_{j} + \hj ,\\
         \xi_{j} = 0 \; \text{on} \; \partial \Omega_{j} - \Gamma_{j} - (\Gamma_{j} + \hj) .
         \end{split}
      \end{equation}
      Define $\vone \defining \sum_{j = 0} ^{I} \grad \xi_{j}\,\ind_{\Omega_{j}}$. Thus, $-\div \vone = \sum_{j = 0} ^{I} \qone\,\ind_{\Omega_{j}}$ and 
      \begin{equation*}
         \vone\cdot\n = \sum_{i = 1} ^{I} \qtwo\,\ind_{\Gammai}
         - \sum_{i = 1} ^{I} \qtwo\,\ind_{\Gammai + \hi}.
      \end{equation*}
      Due to the Poincar\'e inequality $c_{1}\,\Vert \vone \Vert_{\Hdiv (\Omega_{1})} \leq  \Vert \qone\Vert_{L^{2}(\Omega_{1})} + \Vert \qtwo \Vert_{L^{2}(\Gamma)}$. Hence, setting $\vtwo \defining \grad \qtwo$ we have
      \begin{multline}\label{Ineq inf sup inequalities chain}
         \mathcal {B}\v (q) =-\int_{\Omega_{1}}\div\vone\,\qone
+\int_{\Omega_{2}}\vtwo\cdot \grad \qtwo
-\int_{\Gammat} (\vone\cdot\n) \,\qtwo \, d S
+\int_{\Gammab} (\vone\cdot\n) \,\qtwo \, d S\\
= \int_{\Omega_{1}}\vert \qone\vert^{2}
+\int_{\Omega_{2}}\vert \grad \qtwo \vert^{2}
+\int_{\Gammat} \vert \qtwo \vert^{2} \, d S
+\int_{\Gammab} \vert \qtwo \vert^{2} \, d S \\
\geq 
\int_{\Omega_{1}}\vert \qone\vert^{2}
+\frac{\kappa}{2}\int_{\Omega_{2}}\vert  \qtwo \vert^{2}
+ \frac{1}{2} \left(\int_{\Omega_{2}}\vert  \qtwo \vert^{2}
+\int_{\Gamma} \vert \qtwo \vert^{2} \, d S\right) \geq
c \Vert \v\Vert_{\scriptscriptstyle \V} \, \Vert q \Vert_{\scriptscriptstyle Q}
      \end{multline}
For $c \defining \min \{c_{1}, \frac{1}{2}, \frac{\kappa}{2}\}$, which gives the inf-sup condition of the operator $\mathcal{B}$.
   \end{proof}
\end{lemma}
\begin{theorem}\label{Th well-posedness}
Suppose that $0\leq \alpha$, $a_{i}(\cdot)\in L^{\infty} (\Omega)$ and 
\begin{equation}\label{Def minimum coefficient coerciveness}
a^{*} \defining \min _{i \, = \, 1, 2} \essinf \{a_{i}(\x): \x\in \Omega_{i}\} 
\end{equation}
If $a^{*}$ is positive then, the mixed variational formulation \eqref{Pblm mixed formulation} (or equivalently, the system \eqref{Pblm weak form epsilon free system}) is well-posed. 
\begin{proof}
Clearly $\mathcal{A}$ is non-negative and $\V$-coercive on $\ker \mathcal{B}$. The operator $\mathcal {B}$ satisfies the inf-sup condition as seen in the preceding lemma. Due to theorem \eqref{Th General Mixed Formulation well posedness} the result follows.
\end{proof}
\end{theorem}
%
%
%
%
%
%
%
%
\section{Scaling the Problem and Convergence Statements}   
\begin{figure}[!]
\caption[1]{Domains Mapping}\label{Fig Domains Mapping}
\centerline{\resizebox{12cm}{7cm}
{\includegraphics{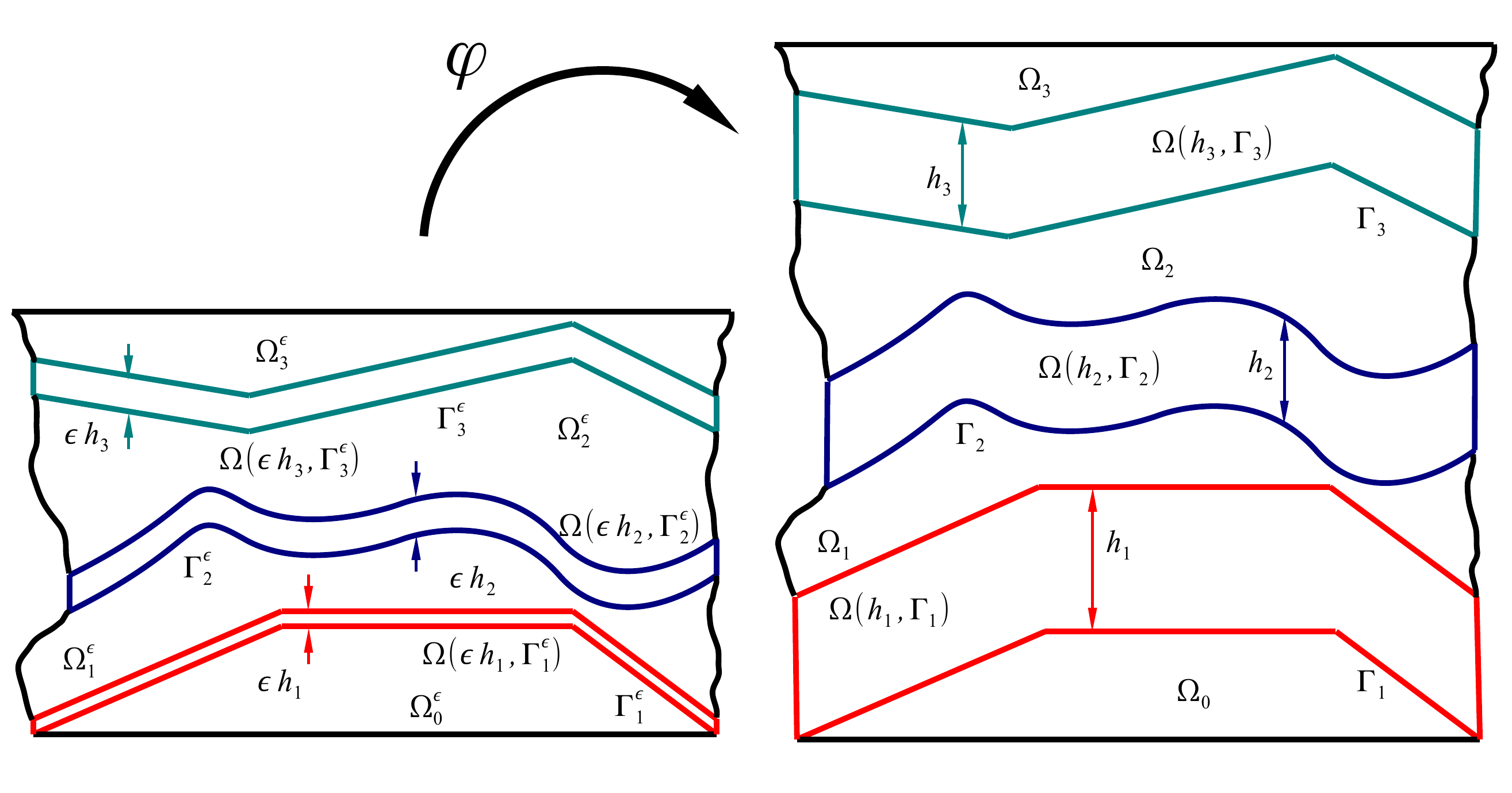} } }
\end{figure}

In order to perform the asymptotic analysis for a the problem \eqref{Pblm weak form epsilon free system} in a medium of thin fractures, the heights and resistance coefficients have to be scaled. We have the following definition (see figure \eqref{Fig Domains Mapping}).
\begin{definition}\label{Def epsilon parametrization fissured medium}
Let $\left\{\left(\zetai, \hi, \Omegai\right):\iindex\right\}$ be a fractured medium of vertical translation generated fissures. For $\epsilon\in (0, 1)$ we define its associated $\epsilon$-scaled fissured  system $\left\{\left(\zetai^{\,\epsilon}, \epsilon\,\hi, \Omegai^{\,\epsilon}\right):\iindex\right\}$ by
\begin{subequations}\label{Eq epsilon fractured medium}
\begin{equation}\label{Eq epsilon fissures surfaces}
\zetai^{\,\epsilon} = \zetai-(1-\epsilon)\sum_{\ell\,=\,0}^{i - 1}h_{\,\ell}\,,\quad 1 \leq i \leq I .
\end{equation}
\begin{equation}\label{Eq epsilon widths}
\left\{\epsilon\,\hi >0 :1 \leq i \leq I \right\} .
\end{equation}
\begin{equation}\label{Eq epsilon shifts of rock matrix}
\Omegaj^{\,\epsilon} \defining \Omegaj-(1-\epsilon)\sum_{\ell\,=\,0}^{j } h_{ \ell}
\,,\quad 0 \leq j \leq I .
\end{equation}
\end{subequations}
The domains $\Omega_{1}^{\epsilon}, \Omega_{2}^{\epsilon}, \Omega^{\epsilon}$ and the surfaces $\Gammat^{\epsilon}, \Gammab^{\epsilon}, \Gamma^{\epsilon}$ are defined as in \eqref{Def rock matrix}, \eqref{Def fissured medium region} respectively. 
\end{definition}
\begin{remark}
   Clearly the systems $\left\{\left(\Gammai^{\epsilon}, \epsilon \, \hi, \Omegai^{\epsilon}\right):\iindex\right\}$ satisfies the conditions of definition \eqref{Def unidirectional translation medium}.  
\end{remark}
%
%
%
%
%
%
\subsection{Isomorphisms of Spaces and Formulation}   
Let $\Omega_{1}^{\epsilon}, \Omega_{2}^{\epsilon}$, $\Omega^{\epsilon}$ and $\Gammat^{\epsilon}, \Gammab^{\epsilon}$, $\Gamma^{\epsilon}$ be the domains and surfaces associated to the family $\{\left(\zetai^{\epsilon}, \epsilon \, \hi, \Omegai^{\epsilon}\right):\iindex\}$ as in definition \eqref{Def epsilon parametrization fissured medium}. Define the spaces
\begin{subequations}\label{Def epsilon spaces of functions}
\begin{equation}\label{Def epsilon space of velocities}
\V^{\epsilon} \defining \{\v\in\Ltwo(\Omega^{\,\epsilon}): \div\vone\in \Ltwo(\Omega_{1}^{\,\epsilon})\,,\vone\cdot\n\,\vert_{\Gamma^{ \epsilon}}\in L^{2}(\Gamma^{ \epsilon}) \} ,
\end{equation}
\begin{equation}\label{Def epsilon space of pressures}
Q^{\epsilon} \defining \{q\in L^{2}(\Omega^{\epsilon}): \grad\qtwo\in \Ltwo(\Omega_{\,2}^{\,\epsilon}) \} .
\end{equation}
We endow the spaces with the norms coming from the natural inner product
\begin{equation}\label{Def epsilon normm space of velocities}
\Vert\, \v\,\Vert_{\scriptscriptstyle \V^{\epsilon}} \defining 
\{\,\Vert \,\v\,\Vert_{\scriptscriptstyle \Ltwo(\Omega^{\epsilon})}^{2}
+\Vert\, \div\vone\,\Vert_{\scriptscriptstyle L^{2}(\Omega_1^{\epsilon})}^{2}
+\Vert\, \vone\cdot\n\,\Vert_{\scriptscriptstyle L^{2}(\Gamma^{\epsilon})}^{2}
\,\}^{1/2}
\end{equation}
\begin{equation}\label{Def epsilon norm of pressures}
\Vert\, q\,\Vert_{\scriptscriptstyle Q^{\epsilon}}
\defining \{\,\Vert\, q\,\Vert_{\scriptscriptstyle L^{2}(\Omega^{\epsilon})}^{2}
+\Vert \,\grad\qtwo\,\Vert_{\scriptscriptstyle L^{2}(\Omega_{\,2}^{\epsilon})}^{2}\,\}^{1/2}
\end{equation}
\end{subequations}
Consider the scaled problem
\begin{subequations}\label{Pblm weak form epsilon system}
\begin{equation*}
\text{Find } \peps\in Q^{\epsilon}, \, \ueps\in \V^{\epsilon} :
\end{equation*}
\begin{multline}\label{Pblm weak form epsilon system 1}
\int_{\Omega_1^{\epsilon}} a_1 \,\ueps \cdot \v \,d\y
+  \epsilon \int_{\Omega_2^{\epsilon}} a_2\, \ueps \cdot \v\,d\y
- \int_{\Omega_1^{\epsilon}} \peps \,\div\v  \,d\y
+ \int_{\Omega_2^{\epsilon}} \grad \peps \cdot\v \,d\y
\\
+\alpha\int_{\Gamma^{\epsilon}} \!\! (\uepsone\cdot\n)\,(\vone\cdot\n ) d S
%
-\int_{\Gammat^{\epsilon}} \!\! \pepstwo (\vone\cdot\n)  d S
+\int_{\Gammab^{\epsilon}} \!\! \pepstwo  (\vone\cdot\n ) d S
= - \!\! \int_{\Omega^{\epsilon}} \!\! \!\g^{\,\epsilon} \cdot \v\,d\y
\end{multline}
\begin{multline}\label{Pblm weak form epsilon system 2}
 \int_{\Omega_1^{\epsilon}}\div\ueps\, q\,d\y
- \int_{\Omega_{2}^{\epsilon}} \ueps \cdot \grad q\,d\y\\
+ \int_{\,\Gammat^{\epsilon}}\left(\uepsone\cdot\n\right) \qtwo\, d S
- \int_{\,\Gammab^{\epsilon}}\left(\uepsone\cdot\n\right) \qtwo\, d S
= \int_{\Omega^{\epsilon}}F^{\,\epsilon}\, q\,d\y
+ \int_{\Gamma^{\epsilon}} f_{\scriptscriptstyle \Gamma^{\epsilon} }^{\,\epsilon}\, \qtwo\,d S
\end{multline}
\begin{equation*}
\text{ for all }  q\in Q^{\epsilon},\,\v\in\V^{\epsilon}
\end{equation*}
\end{subequations}
Clearly, the problem \eqref{Pblm weak form epsilon system} is well-posed since it verifies all the hypothesis of theorem \eqref{Th well-posedness}. In order to analyze the asymptotic behavior of the solution $\left(\ueps, \peps\right)$ as $\epsilon\downarrow 0$ the geometry of the $\epsilon$-domains must be mapped to a common domain of reference. 
%
%
%
%
%
%
\subsection{The $\epsilon$-Problems in a Reference Domain }   
We introduce the change of variable (see figure \eqref{Fig Domains Mapping}) $ \boldsymbol{ \varphi } :\Omega^{\,\epsilon}\rightarrow \Omega$ defined by
%
%
%
%
\begin{multline}\label{Def concise change of variable}
\boldsymbol{\varphi }(\y) \defining \sum_{j = 0}^{I}\left(\widetilde{ \y }, 
y_{\,3}+(1-\epsilon)\sum_{\ell\,=\,0}^{j} h_{\,\ell} \right) \ind_{\Omegaj^{\,\epsilon} }(\y)\\
+ \sum_{i = 0}^{I}\left( \widetilde{ \y }, 
\frac{1}{\epsilon}\left(y_{\,3}-\zetai^{\epsilon}(\widetilde{\y})\right) + \zetai^{\epsilon}(\widetilde{\y})
+(1-\epsilon)\sum_{\ell\,=\,0}^{i - 1} h_{\,\ell} \right) \ind_{\Omega(\epsilon\hi, \Gammai^{\epsilon} ) }(\y)
\end{multline}
%
%
%
Defining $(\xthilde, z) \defining \boldsymbol {\varphi } (\y)$ the gradients are related as follows
\begin{equation}\label{Eq concise gradient structure}
   \grad_{\! \y} = \left\{\begin{array}{c}
            \gradth_{\!\x}\\[7pt]
            \partial_{ z }
            \end{array}\right\}\ind_{\Omega_1^{\epsilon} } + 
   \sum_{i} \left[\begin{array}{cc} 
               I & (1 - \dfrac{1}{\epsilon} ) \, \gradth_{\! \x} \zetai (\xthilde)\\[7pt]
               0 & \dfrac{1}{\epsilon}
            \end{array}\right] \left\{\begin{array}{c}
            \gradth_{\!\x}\\[7pt]
            \partial_{ z }
            \end{array}\right\} \ind_{\Omega(\epsilon \hi, \Gammai^{\epsilon} ) }
\end{equation}
%
%
%
Here, it is understood that $I$ is the identity matrix in $\in \R^{2\times 2}$. We write $\zetai$ instead of $\zetai^{\epsilon}$ for the sake of simplicity recalling that both surfaces differ only by a constant of vertical translation.
\begin{theorem}\label{Th change of variable isomorphism}
Let $\boldsymbol { \varphi }: \Omega^{\,\epsilon}\rightarrow\Omega$ be the change of variable defined in equation \eqref{Def concise change of variable}. Then, the maps defined $\Phi_{1}: \V\rightarrow \V^{\epsilon}$, $\Phi_{\,2}: Q\rightarrow Q^{\epsilon}$ defined respectively by $\left(\Phi_{1}\v\right)(\y) \defining \v\left(\boldsymbol { \varphi } (\y)\right)$ and $\left(\Phi_{\,2}\,q\right)(\y) \defining q\left(\boldsymbol { \varphi } (\y)\right)$ are isomorphisms.
\begin{proof}
First notice for $\v\in \V$ and $q\in Q$ the functions $\Phi_{1}\v$ and $\Phi_{\,2}\,q$ are defined on $\Omega^{\epsilon}$. Moreover, for $\ell = 1, 2$ the restriction of the change of variable is a bijection i.e. $\boldsymbol { \varphi }:\Omega^{\,\epsilon}_{\,\ell}\rightarrow\Omega_{\,\ell}$ is a bijection. Therefore $\v(\cdot)\in \Ltwo(\Omega_{\,\ell})$ if and only if $\v(\boldsymbol { \varphi }(\cdot))\in\Ltwo(\Omega_{\ell}^{\epsilon})$ and $q (\cdot) \in L^{2}(\Omega_{\,\ell})$ if and only if $q(\boldsymbol { \varphi }(\cdot))\in L^{2}(\Omega_{\,\ell}^{\epsilon})$. Even more, $\boldsymbol { \varphi }:\Gammai^{\epsilon}\rightarrow \Gammai$ and $\boldsymbol { \varphi }:\Gammai^{\epsilon} + \epsilon \hi\rightarrow \Gammai + \hi$ are bijective rigid translations. Therefore, the isomorphisms $L^{2}(\Gammai^{\epsilon})\simeq L^{2}(\Gammai)$, $L^{2}(\Gammai^{\epsilon}+\epsilon\,\hi)\simeq L^{2}(\Gammai+\hi)$ follow for all $\iindex$.

For the isomorphism $\Phi_{1}$ take $\v \in \V$ which is equivalent to $\v (\y)\in \Ltwo(\Omega)$ and $\grad_{\y}\cdot\v(\y)\in L^{2}(\Omega_{1})$. Due to the previous discussion these two conditions are equivalent to $\v(\boldsymbol { \varphi }(\y)) \in \Ltwo(\Omega^{\epsilon})$ and $\grad_{\y}\cdot\v(\boldsymbol { \varphi }(\y)) = \grad_{\y}\cdot \v(\x)\in L^{2}(\Omega_{1}^{\epsilon})$. However, equation \eqref{Def concise change of variable} yields $\grad_{\y}\cdot\v(\boldsymbol{\varphi}(\y) ) = \grad_{\y}\cdot \v(\x) = \grad_{\x}\cdot\v(\x)$ whenever $\x\in \Omega_{1}$; i.e. $\grad_{\y}\cdot\v(\y)\in L^{2}(\Omega^{ \epsilon}_{1})$ if and only if $\grad_{\x}\cdot\v(\x) = \grad_{\x}\cdot\v(\boldsymbol {\varphi} (\y))\in L^{2}(\Omega_{ 1})$ as desired.

For the map $\Phi_{\,2}$, the $L^{2}$-integrability condition between spaces $Q$ and $Q^{\epsilon}$ is shown using the same arguments of the first paragraph. It remains to show the $L^{2}$-integrability condition on the gradient. First observe that the last row in the matrix equation \eqref{Eq concise gradient structure} implies that $\frac{\partial}{\partial y_{3}}\, q(\y)\in L^{2}(\Omega_{2}^{\epsilon})$ if and only if $\frac{\partial}{\partial z}\, q(\x)\in L^{2}(\Omega_{2})$. Second, for the derivatives in the first two directions equation \eqref{Eq concise gradient structure} yields
\begin{equation*}
\frac{\partial}{\partial\, y_{\ell}}\, q(\y) = \frac{\partial}{\partial\, x_{\ell}}\, q(\x)+
\left(1-\frac{1}{\epsilon}\right)\frac{\partial}{\partial\, x_{\ell}}\,\zetai(\x)\;\frac{\partial}{\partial z} \, q(\x)\,,\quad \ell = 1, 2 .
\end{equation*}
Recalling the gradient of $\zetai$ is bounded we conclude $\frac{\partial}{\partial y_{\ell}}\, q(\y)\in L^{2}(\Omega_{2}^{\epsilon})$ if an only if $\frac{\partial}{\partial x_{\ell}}\, q(\x)\in L^{2}(\Omega_{2})$ for $\ell = 1, 2$. Since $\frac{\partial}{\partial  z}\, q(\x)\in L^{2}(\Omega_{2})$ is immediate the proof is complete.
\end{proof}
\end{theorem}
We are to apply the change of variable  $\boldsymbol { \varphi }: \Omega^{\,\epsilon}\rightarrow \Omega$ in the problem \eqref{Pblm weak form epsilon system}, to this end, it is more convenient to write the system in terms of the quantities and directions which yield estimates agreeable with the asymptotic analysis. Hence, recalling the definition of the upwards normal vector \eqref{Def Standarazied upwards normal vector} the following relationships hold
\begin{subequations}\label{Eq normal and tangential velocities}
\begin{equation}\label{Eq normal velocity}
\vert(-\widetilde{\grad} \zetai, 1)\vert\,\v\cdot\ni = -\vtan\cdot\widetilde{\grad}\zetai + \vnorm ,
\end{equation}
\begin{equation}\label{Eq vector orthogonal to n}
(\widetilde{\v}, \widetilde{\v}\cdot\gradth\zetai)\cdot\ni = 0
\quad \text{in}\;\, \Omega(\hi, \Gammai) .
\end{equation}
\end{subequations}
Applying the change of variable \eqref{Def concise change of variable} to the problem  \eqref{Pblm weak form epsilon system} and combining with the relation \eqref{Eq normal velocity} we get the following variational statement:
\begin{subequations}\label{Pblm associated weak form epsilon system}
\begin{equation*}
\text{ Find } \peps\in Q, \, \ueps\in \V :
\end{equation*}
\begin{multline}\label{Pblm associated weak form epsilon system 1}
\int_{\Omega_1} a_1 \,\ueps \cdot \v 
+  \epsilon^{\,2} \int_{\Omega_2} a_2\, \ueps \cdot \v
- \int_{\Omega_1} \peps \,\div\v  
\\
+ \sum_{i} \int_{\Omega(\hi, \Gammai)} \epsilon
\left(\widetilde{\grad} \peps+\partial_{z}\,\peps\,\widetilde{\grad} \zetai\right) \cdot\vtan 
+  \vert (-\widetilde{\grad} \zetai, 1 ) \vert \, \partial_{z}\,\peps\,(\v\cdot\ni)
\\
-\int_{\,\Gammat}\pepstwo \left(\vone\cdot\n\right)\, d S
+\int_{\,\Gammab}\pepstwo \left(\vone\cdot\n\right)\, d S \\
+\alpha\int_{\Gamma}\left(\uepsone\cdot\n\right)\left(\vone\cdot\n\right) d S
= - \int_{\Omega_{1}}\g^{\,\epsilon} \cdot \v
- \epsilon \int_{\Omega_{2}}\g^{\,\epsilon} \cdot \v
\end{multline}
\begin{multline}\label{Pblm associated weak form epsilon system 2}
 \int_{\Omega_1}\div\ueps\, q 
- \sum_{ i }\int_{\Omega(\hi, \Gammai)}
\epsilon \, \utaneps \cdot\left(\gradth q+\partial_{z}\,q\,\widetilde{\grad} \zetai\right) 
+ 
\vert (-\widetilde{\grad} \zetai, 1 ) \vert \, (\uepstwo\cdot\ni) \,\partial_{z}\,q \\
+ \int_{\,\Gammat}\left(\uepsone\cdot\n\right) \qtwo\, d S
- \int_{\,\Gammab} (\uepsone\cdot\n )\, \qtwo\, d S 
= \int_{\Omega_{1}}F^{\,\epsilon, 1}\, q
+ \epsilon \int_{\Omega_{2}}F^{\,\epsilon, 2}\, q
+ \int_{\Gamma} f_{\scriptscriptstyle \Gamma}^{\,\epsilon}\, \qtwo\,d S
\end{multline}
\begin{equation*}
\text{ for all }  q\in Q,\,\v\in\V
\end{equation*}
\end{subequations}
Finally, due to the theorem \eqref{Th change of variable isomorphism} on isomorphisms of function spaces we conclude that the problems \eqref{Pblm associated weak form epsilon system} and \eqref{Pblm weak form epsilon system} are equivalent.
%
%
%
%
\subsubsection{The Strong Rescaled Problem}\label{Sec strong epsilon problem}
The solution of the problem \eqref{Pblm associated weak form epsilon system} is the weak solution of the following system of equations
\begin{subequations}\label{Eq epsilon problem strong form}
\begin{equation}\label{Eq epsilon Darcy Omega 1}
a_1 \uepsone+\grad \,\pepsone+\g=0  \text{ and } 
\end{equation}
\begin{equation}\label{Eq epsilon divergence Omega 1}
\div\uepsone=f^{\,\epsilon,\,1} \text{ in }
\Omega_{1} .
\end{equation}
\begin{equation}\label{Eq epsilon drained rock condition rock matrix}
\pepsone = 0 \text{ on } \partial \Omega_{1}-\Gamma .
\end{equation}
\begin{equation}\label{Eq epsilon normal stress balance interface condition}
\pepsone-\pepstwo = 
\alpha \, \uone\cdot\n \; \ind_{\Gammab}
-\alpha \, \uone\cdot\n \; \ind_{\Gammat}
\quad \text{and}
\end{equation}
\begin{equation}\label{Eq epsilon normal flux balance interface condition}
\left(\uepsone -
\uepstwo\right)\cdot\n \; \ind_{\Gammat}
-\left(\uepsone -
\uepstwo\right)\cdot\n \; \ind_{\Gammab} = f^{\epsilon}_{\scriptscriptstyle \Gamma}
\quad \text{on} \;\; \Gamma .
\end{equation}
\begin{equation}\label{Eq epsilon tangential Darcy Omega 2}
\sum_{i}\left[\epsilon\, a_2 \,\utaneps+\widetilde{\grad}\pepstwo
+  (1-\frac{1}{\epsilon} )\, \partial_{\,z}\,\pepstwo\,\widetilde{\grad}\zetai
+\widetilde{\g}^{\,\epsilon} \right]\ind_{\Omega(\hi, \Gammai)} = 0 \, ,
\end{equation}
\begin{equation}\label{Eq epsilon normal Darcy Omega 2}
\epsilon^{\,2} \,a_2 \,\unormeps+\partial_{z}\,\pepstwo+\epsilon\,g^{\,\epsilon}_{3}
= 0 \quad \text{and}
\end{equation}
\begin{multline}\label{Eq epsilon divergence Omega 2}
\sum_{i}\div \left(\epsilon \, \utaneps, \epsilon \, \utaneps\cdot \gradth \zetai + \vert (-\gradth \zetai, 1)\vert (\uepstwo\cdot \ni) \right)
 \ind_{\Omega(\hi, \Gammai) }\\
 =\epsilon\,F^{\,\epsilon,\,2}
 \text{ in } \Omega_{2} .
\end{multline}
\begin{equation}\label{Eq epsilon neumann condition on Omega 2}
\utaneps\cdot\outerth^{\,(i)}=0 \text{ on } \partial\,\Omega(\hi, \Gammai)-\Gamma
\quad \text{for all}\;\iindex .
\end{equation}
\end{subequations}
As before equations \eqref{Eq epsilon normal stress balance interface condition}, \eqref{Eq epsilon normal flux balance interface condition} have the separation of cases $\ind_{\Gammab}, \ind_{\Gammat}$ in order to be consistent with the upwards normal vector $\n$. However, the equations \eqref{Eq epsilon normal flux balance interface condition} and \eqref{Eq epsilon neumann condition on Omega 2} need further clarification. We start fixing an index $i\in \{1, \ldots, I\}$ of the sum in the equation \eqref{Pblm associated weak form epsilon system 2}; reordering and integrating by parts yield
\begin{multline*}
-  \int_{\Omega(\hi, \Gammai)}\epsilon \, \utaneps\cdot(\gradth q + \partial_{z}\,q\, \gradth \zetai)
+ \vert (-\gradth \zetai, 1)\vert (\uepstwo\cdot \ni) \partial_{z} q\\
=  - \int_{\Omega(\hi, \Gammai)}\left(\epsilon \, \utaneps, \epsilon \, \utaneps\cdot \gradth \zetai + \vert (-\gradth \zetai, 1)\vert (\uepstwo\cdot \ni) \right)\cdot\grad q \\
=  \int_{\Omega(\hi, \Gammai)}\div \left(\epsilon \, \utaneps, \epsilon \, \utaneps\cdot \gradth \zetai + \vert (-\gradth \zetai, 1)\vert (\uepstwo\cdot \ni) \right) q\\
-  \int_{\partial \Omega(\hi, \Gammai)}q \left(\epsilon \, \utaneps, \epsilon \, \utaneps\cdot \gradth \zetai + \vert (-\gradth \zetai, 1)\vert (\uepstwo\cdot \ni) \right)\cdot\outer^{(i)} \, dS .
\end{multline*}
Where $\outer^{(i)}$ is the outwards pointing unit normal field of the boundary $\partial \Omega(\hi, \Gammai)$. We focus on the boundary term 
\begin{multline*}
\int_{\partial \Omega(\hi, \Gammai)}q \left(\epsilon \, \utaneps, \epsilon \, \utaneps\cdot \gradth \zetai + \vert (-\gradth \zetai, 1)\vert (\uepstwo\cdot \ni) \right)\cdot\outer^{(i)} \, dS\\
=  \int_{\partial \Omega(\hi, \Gammai) - (\Gammai \cup \hi + \Gammai) }q \left(\epsilon \, \utaneps, \epsilon \, \utaneps\cdot \gradth \zetai + \vert (-\gradth \zetai, 1)\vert (\uepstwo\cdot \ni) \right)\cdot\outer^{(i)} \, dS\\
+ \sum_{\ell \, = \, 0, 1} \int_{\ell \, \hi + \Gammai }q \left(\epsilon \, \utaneps, \epsilon \, \utaneps\cdot \gradth \zetai + \vert (-\gradth \zetai, 1)\vert (\uepstwo\cdot \ni) \right)\cdot\outer^{(i)} \, dS .
\end{multline*}
The equality $\outer^{(i)}\cdot\kversor = 0$ holds on the portion of the vertical wall $\partial \Omega(\hi, \Gammai) - (\Gammai \cup \hi + \Gammai)$ i.e. the equation \eqref{Eq epsilon neumann condition on Omega 2} follows. For the remaining pieces of the boundary recall $\ni = \outer^{(i)}$ on $\hi + \Gammai$ and $\ni = - \outer^{(i)}$ on $ \Gammai$; together with the equation \eqref{Def Standarazied upwards normal vector}, we get
\begin{multline*}
- \int_{\ell \, \hi + \Gammai }q \left(\epsilon \, \utaneps, \epsilon \, \utaneps\cdot \gradth \zetai + \vert (-\gradth \zetai, 1)\vert (\uepstwo\cdot \ni) \right)\cdot\outer^{(i)} \, dS\\
= (-1)^{\ell }\int_{\ell \, \hi + \Gammai } \!\! q \left(\epsilon \, \utaneps, \epsilon \, \utaneps\cdot \gradth \zetai + \vert (-\gradth \zetai, 1)\vert (\uepstwo\cdot \ni) \right)\cdot\frac{(-\gradth \zetai, 1)}{\vert (-\gradth \zetai, 1)\vert} \, dS \\
= (-1)^{\ell }\int_{\ell \, \hi + \Gammai }  q  \, (\uepstwo\cdot \ni) \,  dS \quad\text{for}\;\ell = 0, 1.
\end{multline*}
Combining this last identity with the interface terms in equation \eqref{Pblm associated weak form epsilon system 2}, the strong normal flux balance condition \eqref{Eq epsilon normal flux balance interface condition} follows.
%
%
%
%
\subsection{A-priori Estimates and Convergence Statements}\label{Sec a priori estimates and convergence}
In order to get a-priori estimates on the norm of the solutions the following hypothesis are assumed
\begin{subequations}\label{Eq boundedness forcing terms}
\begin{equation}\label{Eq boundedness forcing divergence terms}
\Vert F^{\,\epsilon}\Vert_{L^{2}(\Omega)} \;\text{is bounded and}\;
F^{1,\epsilon}\overset{w}\rightharpoonup F^{\,1}\;\text{in}\;L^{2}(\Omega_{1}) , 
\end{equation}
\begin{equation}\label{Eq boundedness forcing gravity terms}
\g^{\epsilon}\overset{w}\rightharpoonup \g\;\text{in}\; \Ltwo(\Omega_{1})\,,\;
\g^{\,2,\epsilon}(\xthilde, \epsilon\,z)\overset{w}\rightharpoonup \g(\xthilde)\;\text{in}\; \Ltwo(\Omega_{2}) ,
\end{equation}
\begin{equation}\label{Eq boundedness interface forcing terms}
\text{and}\; f^{\,\epsilon}_{\scriptscriptstyle \Gamma}\overset{w}\rightharpoonup f_{\scriptscriptstyle \Gamma}\;\text{in}\;L^{2}(\Gamma) .
\end{equation}
\end{subequations}
Now test equation \eqref{Pblm associated weak form epsilon system 1} with $\ueps$ and equation \eqref{Pblm associated weak form epsilon system 2} with $\peps$ add them together and get
\begin{multline}\label{Eq first estimate}
a^{*} \left(\Vert\,\uepsone\Vert_{ 0, \Omega_1}^{ 2}
+\Vert\,
\epsilon\,\uepstwo\Vert_{0, \Omega_2}^{ 2}\right)
+\alpha \left\Vert\,
\uepsone\cdot\n\right\Vert_{ L^{2}(\Gamma)}^{ 2}\\
=
\int_{\Omega_{1}}F^{\,\epsilon,\,1}\,\peps 
+\epsilon\int_{\Omega_{ 2}}F^{\,\epsilon,\,2}\,\peps 
+\int_{\Gamma}f^{\,\epsilon}_{\scriptscriptstyle \Gamma}\,\pepstwo\, dS 
-\int_{\Omega_1}\g^{\,1}\cdot\ueps 
-\int_{\Omega_2}\g^{\,2}\cdot \epsilon\,\ueps 
\\
\leq C \left(\Vert F^{\,\epsilon}\Vert_{ 0,\,\Omega}
+\Vert f^{\,\epsilon}_{\scriptscriptstyle \Gamma}\Vert_{ 0,\,\Gamma}\right)\,\Vert \peps\Vert_{Q}\,
%
+\Vert \g^{\,\epsilon}\Vert_{ 0,\,\Omega}
\left(\Vert \uepsone \Vert_{ 0,\,\Omega_1}
+\Vert \epsilon\,\uepstwo \Vert_{ 0,\,\Omega_{ 2}}\right) .
\end{multline}
Here, the constant $C>0$ is independent from $\epsilon>0$. Next the term $\Vert \peps \Vert_{\scriptscriptstyle Q}$ must be bounded in terms of the flux $\uepsone \ind_{\Omega_{1}} + \epsilon \, \uepstwo\ind_{\Omega_{2}}$ and the forcing terms. Due to the equation \eqref{Eq epsilon normal Darcy Omega 2} we have
\begin{subequations}\label{Eq boundedness pressure gradient omega 2}
\begin{equation}\label{Eq boundedness normal pressure gradient omega 2}
\Vert \frac{1}{\epsilon}\,\partial_{z}\, \pepstwo \Vert_{ 0, \Omega_{2}} \leq \epsilon \,\Vert a_{2}\Vert_{\scriptscriptstyle L^{\infty}(\Omega_{2})}
\Vert \unormeps \Vert_{0, \Omega_{2}}+\Vert g^{\epsilon}_{\,3} \Vert_{0, \Omega_{2}} .
\end{equation}
Combined with equation \eqref{Eq epsilon tangential Darcy Omega 2} yields
\begin{equation}\label{Eq boundedness tangential pressure gradient omega 2}
\Vert \gradth \pepstwo \Vert_{ 0, \Omega_{2}} \leq C \left(
\Vert a_{2}\Vert_{\scriptscriptstyle L^{\infty}(\Omega_{2})}
\Vert \epsilon\, \uepstwo \Vert_{0,\,\Omega_{2}}
+\Vert \g^{\,\epsilon} \Vert_{0,\,\Omega_{2}}\right) .
\end{equation}
\end{subequations}
For $C>0$ an adequate constant. Thus
\begin{equation}\label{Eq boundedness pressure full gradient omega 2}
\Vert \grad \pepstwo \Vert_{ 0, \Omega_{2}} \leq  
C \left(\Vert a_{2}\Vert_{\scriptscriptstyle L^{\infty}(\Omega_{2})}
\Vert \epsilon\,\uepstwo \Vert_{0,\,\Omega_{2}}
+\Vert \g \Vert_{0,\,\Omega_{2}}\right).
\end{equation}
With $C>0$ a constant independent from $\epsilon>0$. Additionally, the equation \eqref{Eq epsilon Darcy Omega 1} yields
\begin{equation}\label{Eq boundedness pressure full gradient omega 1}
\Vert \grad \pepstwo \Vert_{ 0, \Omega_{1}} \leq  \Vert a_{1}\Vert_{\scriptscriptstyle L^{\infty}(\Omega_{2})}
\Vert  \uepstwo \Vert_{0, \Omega_{1}} +\Vert \g \Vert_{0, \Omega_{1}} .
\end{equation}
The boundary condition \eqref{Eq epsilon drained rock condition rock matrix} together with Poincar\'e inequality give the control $\Vert  \pepsone \Vert_{H^{1}(\Omega_{1})} \leq C \Vert \grad \peps \Vert_{0, \Omega_{1}}$. On the other hand, the inequality \eqref{Ineq control by gradient and trace} implies $\Vert \peps \Vert_{1, \Omega_{2}} \leq C (\Vert \peps \Vert_{0, \Gamma} + \Vert \peps \Vert_{1, \Omega_{2}}) $; combined with the normal stress balance conditions \eqref{Eq epsilon normal stress balance interface condition} we conclude: 
\begin{equation}\label{Ineq presssure norm controled by pressure gradient on the solution}
\Vert \peps \Vert_{\scriptscriptstyle Q} \leq \Vert \peps \Vert_{1, \Omega} \leq C\,\Vert \grad \peps \Vert_{0, \Omega} .
\end{equation}
And $C>0$ is independent from $\epsilon>0$. Finally, a combination of inequalities \eqref{Ineq presssure norm controled by pressure gradient on the solution}, \eqref{Eq boundedness pressure full gradient omega 1} and \eqref{Eq boundedness pressure full gradient omega 2} imply that the left hand side of inequality \eqref{Eq first estimate} is bounded. 
\begin{remark}\label{Rem different approach a priori estimates}
The previous estimate on $\Vert  \pepstwo \Vert_{H^{1}(\Omega_{2})}$ could have been attained without requiring the drained condition \eqref{Eq epsilon drained rock condition rock matrix} on the whole matrix rock region external boundary. It was enough to set the drained condition on a subset of positive measure contained in $\partial \Omegaj - \Gamma$ for $j$ fixed to have control on $\Vert  \pepsone \Vert_{1, \Omegaj}$ by $\Vert  \grad \pepsone \Vert_{0, \Omegaj}$. Combining this fact with the normal stress balance conditions \eqref{Eq epsilon normal stress balance interface condition}, an inequality of the type \eqref{Ineq presssure norm controled by pressure gradient on the solution} can be deduced for the union of adjacent domains $\Omega(h_{j}, \Gamma_{j}) \cup \Omegaj \cup \Omega(h_{j + 1}, \Gamma_{j + 1})$ and continue the process until the whole domain $\Omega$ is covered and the global inequality \eqref{Ineq presssure norm controled by pressure gradient on the solution} is obtained.
\end{remark}
Due to the observations above we conclude that the following sequences are bounded
\begin{subequations}\label{Eq boundedness of epsilon solutions}
\begin{equation}\label{Eq boundedness epsilon velocities}
\Vert \uepsone \Vert_{ 0, \Omega_{1}} \,,\;
\Vert \epsilon\,\uepstwo \Vert_{0,\,\Omega_{2}}\,,\;
\sqrt{\alpha} \, \Vert \uepsone\cdot\n \Vert_{L^{2}(\Gamma)}
\end{equation}
\begin{equation}\label{Eq boundedness of epsilon pressures}
\Vert \pepsone \Vert_{ H^{1}(\Omega_{1})} \,,\;\Vert \pepstwo \Vert_{H^{1}(\Omega_{2})}\,,\;
\Vert \frac{1}{\epsilon}\,\partial_{z}\,\peps \Vert_{0,\,\Omega_{2}}
\,,\;
\Vert \div \uepsone \Vert_{L^{2}(\Omega_{1})} .
\end{equation}
\end{subequations}
\begin{remark}\label{Rem Divergence Behavior}
The change of variable $\boldsymbol { \varphi }$ modifies the structure of the divergence on the domains $\Omega(\hi, \Gammai)$ for all $1\leq i\leq I$, therefore it can only be claimed that the linear combination $\epsilon\,\widetilde{\grad}\cdot\utaneps
+\epsilon\, (1-\frac{1}{\epsilon} )\partial_{\,z} (\widetilde{\grad}\zetai\cdot\utaneps )
+ \partial_{z}\unormeps$ is bounded in $L^{2}(\Omega(\hi, \Gammai))$.
\end{remark}
%
%
%
%
\subsection{Weak Limits}\label{Sec weak convergence}
The previous section state bounds independent from $\epsilon>0$ for $[\uepsone, \epsilon\,\uepstwo]\in \V$ and $\peps = [\pepsone, \pepstwo]\in H^{1}(\Omega_{1})\times H^{1}(\Omega_{2})$, consequently in $Q$. Then, there must exist $\u\in \V$, $p\in Q$, $\eta\in L^{2}(\Omega_{2})$ and a subsequence, from now on denoted the same, such that
\begin{subequations}\label{Eq convergence of epsilon solutions}
\begin{equation}\label{Eq convergence of epsilon pressures}
\peps \overset{w}\rightharpoonup p\;\,\text{in}\,\;Q\,\; \text{and strongly in}\;L^{2}(\Omega) ,
\end{equation}
\begin{equation}\label{Eq convergence of epsilon velocities 1}
\uepsone \overset{w}\rightharpoonup \uone\;\,\text{in}\,\;\Ltwo(\Omega_{1})\,\;
\text{and}\,\;\div\uepsone \overset{w}\rightharpoonup \div\uone\;\,\text{in}\,\;L^{2}(\Omega_{1}) ,
\end{equation}
\begin{equation}\label{Eq convergence of epsilon velocities on interface}
\sqrt{\alpha}\,\uepsone\cdot\n \overset{w}\rightharpoonup \sqrt{\alpha}\,\uone\cdot\n\;\,\text{in}\,\;L^{2}(\Gamma) ,
\end{equation}
\begin{equation}\label{Eq convergence of epsilon velocities 2}
\epsilon\, \uepsone \overset{w}\rightharpoonup \utwo\;\,\text{in}\,\;\Ltwo(\Omega_{2}) ,
\end{equation}
\begin{equation}\label{Eq convergence of higher order pressure}
\frac{1}{\epsilon}\,\partial_{z}\pepstwo \overset{w}\rightharpoonup \eta\;\,\text{in}\,\;L^{2}(\Omega_{2})\,\; \text {and}\,\;
\partial_{z}\pepstwo \rightarrow 0\;\,\text{strongly in}\,\;L^{2}(\Omega_{2}) .
\end{equation}
\end{subequations}
Choose $\phi\in C_{0}^{\infty}(\Omega(\hi, \Gammai))$ arbitrary, test the equation \eqref{Pblm associated weak form epsilon system 2} with $q \defining \epsilon \phi$ and let $\epsilon\downarrow 0$. Recalling \eqref{Eq convergence of epsilon velocities 2} this gives
\begin{multline*}
0 = \lim_{\epsilon\downarrow 0} \int_{\Omega(\hi, \Gammai)} \vert(-\gradth \zetai, 1) \vert \, (\epsilon\,\uepstwo\cdot\ni)\,\partial_{z}\phi 
= \int_{\Omega_{2}} \vert(-\gradth \zetai, 1) \vert \, (\utwo\cdot\ni)\,\partial_{z}\phi\\
= -\left\langle \partial_{z} \vert(-\gradth \zetai, 1) \vert \, (\utwo\cdot\ni), \phi\right\rangle _{D\,'(\Omega(\hi, \Gammai)), D(\Omega(\hi, \Gammai))} .
\end{multline*}
Since $(-\gradth \zetai, 1)$ does not depend on the vertical variable $z$ and it is the non-zero vector almost everywhere we conclude $\partial_{z}(\utwo\cdot\ni) = 0$ i.e. the component of the velocity normal to the surface $\Gammai$ is independent from $z$ in $\Omega(\hi, \Gammai)$ for all $\iindex$. Now choose $q\in Q$ arbitrary, test \eqref{Pblm associated weak form epsilon system 2} with $\epsilon\, q$ and let $\epsilon \downarrow 0$ to get
\begin{multline*}
0 = \sum_{i}\int_{\Omega(\hi, \Gammai)} \vert(-\gradth \zetai, 1) \vert 
\, (\utwo\cdot\ni)\,\partial_{z} q \; d\x\\
= \sum_{i}\int_{\Gi}\int_{\zetai(\xthilde)}^{\,\zetai(\xthilde)+\hi}
\vert(-\gradth \zetai, 1)\vert \, (\utwo\cdot\ni)\,\partial_{z} q\; dz\,d\xthilde\\
%
%
= \sum_{i}  \int_{\Gi} \vert(-\gradth \zetai, 1) \vert(\utwo\cdot\ni)
\left[ q (\xthilde, \zetai(\xthilde)+\hi) - q(\xthilde, \zetai(\xthilde))\right]\,d\xthilde .
\end{multline*}
The above holds for all $q\in Q$, in particular choosing $q(\xthilde, \zetai(\xthilde)) = \phi(\xthilde)$ for $\phi\in C_{0}^{\infty}(\Gi)$ arbitrary and $q(\xthilde, \zetai(\xthilde)+\hi) = 0$ the statement transforms in
\begin{equation*}
%
\int_{\Gi} \vert(-\gradth \zetai, 1) \vert\, 
(\utwo\cdot\ni)(\xthilde,\,\zetai(\xthilde))\,\,\phi\,(\xthilde,\,\zetai(\xthilde))\,\,d\xthilde
%
\quad\forall\,\phi\in C_{0}^{\infty}(\Gi) .
\end{equation*}
Therefore $\vert(-\widetilde{\grad} \zetai, 1)\vert
\left(\utwo\cdot\n\right)$ must be null and since $\vert(-\widetilde{\grad} \zetai, 1 ) \vert$ is non-zero almost everywhere we conclude
\begin{equation}\label{Eq null normal velocity}
\utwo\cdot\ni = 0 \;\,\text{in} \;\,\Omega(\hi, \Gammai)\; \text{for each} \; \iindex .
\end{equation}
The later implies that the Cartesian coordinates of $\utwo$ satisfy the following relation
\begin{equation}\label{Eq structure of limit velocity}
\utwo = \left\{
\begin{array}{c}
\utan \\
\unorm
\end{array}\right\}=
\left\{
\begin{array}{c}
\utan \\
\utan\cdot\widetilde{\grad}\zetai
\end{array}\right\}\;\,\text{in} \;\,\Omega(\hi, \Gammai)\,,\;\iindex .
\end{equation}
Now fix $i\in \{1,\ldots, I\}$ and take a function $\vtang\in(C_{0}^{\infty}(\Omega(\hi, \Gammai))^{2}$. Recalling \eqref{Eq local and global velocities} define $\vtan\defining \MTtang_{i} \vtang$ and $\vnorm\defining \MNtang \vtang$. Then, the function $\displaystyle \vtwo\defining \frac{1}{\epsilon}\,(\vtan,\,\vnorm)$ has the structure \eqref{Eq structure of limit velocity} or equivalently $\vtwo\cdot\n = 0$ inside $\Omega(\hi, \Gammai)$. Define $\v$ as the trivial extension of $\vtwo$ to the whole domain $\Omega$, therefore $\v\in \V$. Test \eqref{Pblm associated weak form epsilon system 1} with $\v$ and let $\epsilon \downarrow 0$, this gives
\begin{equation*}
\int_{\Omega(\hi, \Gammai)}a_2\,(\x)\,\utwo\cdot(\vtan,\,\vnorm) 
+\int_{\Omega(\hi, \Gammai)}\widetilde{\grad}\,\ptwo
\cdot\vtan 
+\int_{\Omega(\hi, \Gammai)}\g\cdot(\vtan,\,\vnorm)  = 0 .
\end{equation*}
Consequently
\begin{equation*}
\int_{\Omega(\hi, \Gammai)}a_2\,(\x)\,\utang\cdot\v_{\tau}
+\int_{\Omega(\hi, \Gammai)}\widetilde{\grad}\,\ptwo
\cdot\MTtang_{i}\v_{\tau}
+\int_{\Omega(\hi, \Gammai)}\g_{\,\tau}\cdot \v_{\tau} = 0 .
\end{equation*}
The equation above holds for all $\v_{\tau} \in (C_{0}^{\infty}(\Omega(\hi, \Gammai))^{2}$ and due to the isomorphism of proposition \eqref{Th isometry change of coordinates} we conclude 
\begin{equation}\label{Eq lower dimensional Darcy-type}
a_2\,(\x)\,\utang+ \left(\MTtang_{i}\right)' \,\widetilde{\grad}\,\ptwo
+\g_{\,\tau} = 0 \,\; \text{in}\; \Omega(\hi, \Gammai) ,\;\,\iindex .
\end{equation}
The equation \eqref{Eq convergence of higher order pressure} implies that $\ptwo$ does not depend on the variable $z$ on $\Omega_{2}$ i.e. $\ptwo = \ptwo(\xthilde)$. Therefore assuming
\begin{equation}\label{Eq dependence of limit data}
a_2 = a_{2}(\xthilde)\,,\quad \widetilde{\g} = \widetilde{\g}(\xthilde)\,\,\;\text{in}\;\,\Omega_{2}
\end{equation}
the equation \eqref{Eq lower dimensional Darcy-type} gives $\utang = \utang(\xthilde)$ i.e. $\utang$ is independent from $z$ in $\Omega_{2}$. Together with the fact $\utwo\cdot\n = 0$ in $\Omega_{2}$ we conclude that the whole vector velocity $\utwo$ is independent from $z$ in $\Omega_{2}$.
\begin{remark}\label{Rem lower dimensional Darcy's law}
   Observe that due to the assumptions for the data \eqref{Eq dependence of limit data} the equation \eqref{Eq lower dimensional Darcy-type} is independent from $z$ becoming a lower-dimensional Darcy-type constitutive law on the stream lines parallel to $\zetai$.
\end{remark}

%
%
%
%
\section{The Limit Problem}\label{Sec limit problem}
Define the subspaces
\begin{subequations}\label{Def limiti space of functions}
\begin{equation}\label{Def limit space of velocities}
\V_{\!0} \defining \left\{\v\in \V: \partial_{z}\vtwo = 0\,,\; \vtwo\cdot\ni = 0 \;\text{in}\;\Omega(\hi, \Gammai) \right\} .
\end{equation}
\begin{equation}\label{Def limit space of pressures}
Q_{0} \defining \left\{q\in Q: \partial_{z}\,q = 0\;\text{in}\; \Omega_{2}\right\} .
\end{equation}
\end{subequations}
\begin{remark}\label{Rem traces of limit pressure}
   Notice that if $q\in Q_{0}$ the fact that $\partial_{z}\,q = 0$ in $\Omega_{2}$ implies $\qtwo\vert_{\Gammai} = \qtwo\vert_{\Gammai+\hi}$ for all $\iindex$. Therefore, the spaces
\begin{equation*}
\begin{split}
\left\{\vone\in \Ltwo(\Omega_{1}):\div\vone\in L^{2}(\Omega_{1}),\,\vone\cdot\ni\in L^{2}(\Gammai)\right\} & \times
\prod_{ i  } L^{2}(\Gammai) \\
%
L^{2}(\Omega_{1})\times\prod_{ i } H^{1}(\Omega(\hi, \Gammai)) &
\end{split}
\end{equation*}
Are isomorphic to \eqref{Def limit space of velocities} and \eqref{Def limit space of pressures} respectively.
\end{remark}
Due to the structure of the space if $\v = [\vone, \vtwo]\in \V_{0}$ then the function $[\vone, \frac{1}{\epsilon}\,\vtwo]$ is also in $\V_{0}$. Using the latter to test \eqref{Pblm associated weak form epsilon system 1} and $q\in Q_{0}$ for testing \eqref{Pblm associated weak form epsilon system 2} we let $\epsilon\downarrow 0$ and conclude that the limits $[\uepsone, \epsilon\,\uepstwo]\rightarrow \u$ and $\peps\rightarrow p$ are a solution of the \emph{limit problem}
\begin{subequations}\label{Pblm limit system}
\begin{equation*}
\text{Find} p\in Q_{0}, \, \u \in \V_{\!0} :
\end{equation*}
\begin{multline}\label{Pblm limit system 1}
\int_{\Omega_1} a_1 \,\u \cdot \v  
- \int_{\Omega_1} p \,\div\v   
+  \int_{\Omega_2} a_2\, \utang \cdot \vtang 
%
+ \int_{\Omega_{2}}
\widetilde{\grad} \peps \cdot\vtan  \\
-\int_{\Gammat} \!\! \ptwo ( \vone\cdot\n )  d S
+\int_{\Gammab} \!\! \ptwo ( \vone\cdot\n )  d S 
+\alpha \!\! \int_{\Gamma} ( \uone\cdot\n) (\vone\cdot\n ) d S 
= - \int_{\Omega_{1}} \!\! \g \cdot \v 
-  \int_{\Omega_{2}} \!\! \g_{\,\tau} \cdot \vtang 
\end{multline}
\begin{multline}\label{Pblm limit system 2}
 \int_{\Omega_1}\div\u\; q 
- \int_{\Omega_{2}}
\utan \cdot \gradth q  
+ \int_{\,\Gammat}\left(\uone\cdot\n\right) \qtwo\, d S
- \int_{\,\Gammab}\left(\uone\cdot\n\right) \qtwo\, d S \\
= \int_{\Omega_{1}}F^{\, 1}\, q 
+ \int_{\Gamma} f_{\scriptscriptstyle \Gamma}\, \qtwo\,d S
\end{multline}
\begin{equation*}
\text{ for all }  q\in Q_{0},\,\v\in\V_{\! 0}
\end{equation*}
\end{subequations}
%
%
%
%
\subsection{Well-Posedness of the Limit Problem}\label{Sec well posedness limit problem}
The problem \eqref{Pblm limit system} is a mixed formulation of the type \eqref{Pblm mixed formulation} with the operators $\mathcal{A}^{0}:\V_{\!0}\rightarrow\V_{\!0}'$ and $\mathcal{B\,}^{0}:
V_{\!0} \rightarrow Q_{0} '$ defined by
\begin{subequations}\label{Def limit bilinear forms}
\begin{equation}\label{Def limit bilinear form velocities}
\mathcal{A}^{0}\v(\w) \defining \int_{\Omega_{1}}\aone\,\v\cdot\w 
+\int_{\Omega_{2}}\atwo\,\vtang\cdot\wtang  +
\alpha \int_{\Gamma}\left(\vone\cdot\n\right)\left(\wone\cdot\n\right)\,d S
\end{equation}
\begin{equation}\label{Def limit bilinear form mixed}
\mathcal{B}^{\,0}\v(q) \defining -\int_{\Omega_{1}}\div\v \; q 
+\int_{\Omega_{2}}\vtan\cdot \gradth q 
-\int_{\Gammat} (\vone\cdot\n) \; \qtwo 
+\int_{\Gammab} (\vone\cdot\n) \; \qtwo 
\end{equation}
\end{subequations}
\begin{theorem}\label{Th limit inf-sup condition}
The operator $\mathcal{B}^{\,0}$ satisfies the inf-sup condition.
\begin{proof}
The proof has the same structure as lemma \eqref{Th inf-sup condition original}, there is only one detail to be examined in the construction of the test functions. Fix $q = [\qone, \qtwo]\in Q_{0}$, construct $\vone$ in the same way it is built in problem \eqref{Pblm Test Function}. On the other hand since $\qtwo \in H^{\,1}\left(\Omega_{\,2}\right)$, $\partial_z\, \qtwo = 0$, define 
\begin{equation*}
\vtwo\defining \sum_{i} (\widetilde{\grad}\,\qtwo, \vtan\cdot\widetilde{\grad}\,\zetai) 
\ind_{\Omega (\hi, \Gammai)}.
\end{equation*}
Then $\vtwo\cdot\ni = 0$ and $\partial_{z} \vtwo = 0$ in $\Omega(\hi, \Gammai)$ for all $\iindex$, i.e. $\vtwo\in \V_{\! 0}$ and $\Vert\,\vtwo\,\Vert_{ 0,\Omega_2}\,\leq C\,\Vert\,\qtwo\,\Vert_{ 1,\Omega_2}$ as desired. Repeating the inequalities presented in \eqref{Ineq inf sup inequalities chain} the proof is complete.
\end{proof}
\end{theorem}
Since the inf-sup condition holds the theorem \eqref{Th well-posedness} applies to the operators \eqref{Def limit bilinear forms} on the spaces $\V_{\! 0}, Q_{0}$ and the limit problem \eqref{Pblm limit system} is well-posed. Due to the uniqueness of the solution of the limit problem it follows that the original sequence converges weakly to the limit $\u\in \V_{\! 0},\, p\in Q_{0}$.
%
%
\subsection{The Strong Form}\label{Sec strong limit problem}
In order to describe the \emph{strong limit problem} corresponding to \eqref{Pblm limit system} two features have to be exploited. First, the structure $\v\cdot\ni = 0$ in $\Omega_{2}$ for all $\v\in \V_{\! 0}$ implying $\vtan = \MTtang \vtang$, for $\MTtang$ the matrix defined in \eqref{Eq local and global velocities}. Second, the independence of the velocities and pressures with respect to $z$ in $\Omega_{2}$. This last property allows to write the integrals over $\Omega(\hi, \Gammai)$ as surface integrals on $\Gammai$. Hence, the system \eqref{Pblm limit system} transforms in
\begin{subequations}\label{Pblm lower dimensional variational statement}
\begin{equation*}
\text { Find } p\in Q_{0}\,,\; \u\in \V_{\! 0} :
\end{equation*}
\begin{multline}\label{Pblm lower dimensional variational statement 1}
\int_{\Omega_1} a_1 \,\u \cdot \v 
- \int_{\Omega_1} p \,\div\v  
+ \sum_{i}\hi \!\! \int_{\Gammai} \!\! (\ni\cdot\kversor\,)(a_2\, \utang + (\MTtang_{i})'\,
\widetilde{\grad} \peps +\g_{\,\tau})\cdot \vtang \,d S\\
-\int_{\,\Gammat}\ptwo \left(\vone\cdot\n\right)\, d S
+\int_{\,\Gammab}\ptwo \left(\vone\cdot\n\right)\, d S 
+\alpha\int_{\Gamma}\left(\uone\cdot\n\right)\left(\vone\cdot\n\right) d S
= - \int_{\Omega_{1}}\g \cdot \v 
\end{multline}
\begin{multline}\label{Pblm lower dimensional variational statement 2}
 \int_{\Omega_1}\div\u \; q 
- \sum_{i}\hi\!\!\int_{\Gammai}\!\! (\ni\cdot\kversor\,)\,
\MTtang_{i}\utang \cdot \gradth \qtwo \,d S\\
+ \int_{\,\Gammat}\left(\uone\cdot\n\right) \qtwo\, d S
- \int_{\,\Gammab}\left(\uone\cdot\n\right) \qtwo\, d S
= \int_{\Omega_{1}}F^{\, 1}\, q 
+ \int_{\Gamma} f_{\scriptscriptstyle \Gamma}\, \qtwo\,d S
\end{multline}
\begin{equation*}
\text{for all}\;q\in Q_{0}\,,\; \v\in \V_{\! 0}
\end{equation*}
\end{subequations}
Integrating by parts the above statement we get the strong lower dimensional problem
\begin{subequations}\label{Pblm limit problem strong form}
\begin{equation}\label{Pblm strong limit Darcy Omega 1}
a_1\,\u+\grad \pone+\g^{1}=0,
\end{equation}
\begin{equation}\label{Pblm strong limit divergence Omega 1}
\div\u=F^{\,1} \text{ in } \Omega_{1} .
\end{equation}
\begin{equation}
\pone = 0 \text{ on } \partial \Omega_{1}-\Gamma .
\end{equation}
\begin{equation}
\utwo\cdot\n=0\,, 
\quad \partial_z\,\ptwo=0\,,
\end{equation}
\begin{equation}\label{Pblm strong limit Darcy Omega 2}
\sum_{i} \left[ a_2\,(s)\,\utwo_{\,\tau}+(\MTtang_{i} )'\,\widetilde{\grad}\ptwo +\g_{\,\tau}^{\,2} (s) 
\right] \ind_{\Gammai}= 0\,,
\end{equation}
\begin{multline}\label{Pblm strong limit divergence Omega 2}
\sum_{i} \left[ (\uone\cdot\ni\vert_{\,\Gammai + \hi}-\uone\cdot\ni\vert_{\,\Gammai}) \right]\ind_{\Gammai}\\
+ \sum_{i}\hi\,  (\ni\cdot\kversor\,)\,\widetilde{\grad}\cdot(\MTtang_{i} \utwo_{\,\tau}) \,  \ind_{\Gammai}
= f_{\scriptscriptstyle\Gamma}\,\;
\text{in}\;\Gamma .
\end{multline}
\begin{equation}\label{Pblm strong normal stress balance}
\pone - \ptwo =
\alpha \, \uone\cdot\n \, \ind_{\Gammab}
-\alpha \, \uone\cdot\n \,  \ind_{\Gammat} .
\end{equation}
\begin{equation}\label{Pblm strong limit flux boundary condition}
\utwo\cdot\outer^{(i)} =  0\,\;\text{on}\;\partial\,\Gi\;\,
\text{for all}\;\iindex .
\end{equation}
\end{subequations}
The statement of equation \eqref{Pblm strong limit Darcy Omega 2} was already shown in \eqref{Eq lower dimensional Darcy-type}, however the statements \eqref{Pblm strong limit divergence Omega 2} and \eqref{Pblm strong limit flux boundary condition} need further discussion.
%
%
\subsection{The Interface Integrals Setting}\label{Sec interface integrals setting}
\begin{definition}\label{Def manifold spaces}
Let $G$, $\zeta$ and $\Gamma$ be as in definition \eqref{Def eligible surface and generated fissure}, define the spaces
\begin{subequations}\label{Def lower dimensional manifold isomorphisms}
\begin{equation}\label{Def L2 manifold}
L^{2}(\Gamma) \defining \{h:\Gamma\rightarrow\R: \int_{\Gamma} h^{2}(s) \,d S <+\infty  \}
\end{equation}
\begin{equation}\label{Def H1 manifold}
H^{1}(\Gamma) \defining \{h  \in L^{2}(\Gamma): \gradth h \in L^{2}(\Gamma)\times L^{2}(\Gamma)\}
\end{equation}
\begin{equation}\label{Def H zero one manifold}
H^{1}_{0}(\Gamma) \defining \{h\in H^{1}(\Gamma): h\vert_{\partial\,G} = 0\}
\end{equation}
\end{subequations}
Here $\gradth$ indicates the gradient with respect to the variables $(x_{1}, x_{2})$ contained in $G$.
\end{definition}
The following isomorphism result is necessary
\begin{theorem}\label{Th manifold isomorphism}
Let $G$, $\zeta$ and $\Gamma$ be as in definition \eqref{Def eligible surface and generated fissure}. Consider the natural embedding $\jmath: \Gamma \rightarrow G$ defined by $\jmath (\xthilde, \zeta(\xthilde)) \defining \xthilde$ and the map
\begin{equation}\label{Def embedding}
\varphi \mapsto \varphi \circ \jmath
\end{equation}
Then
\begin{enumerate}[(i)]
\item
The embedding \eqref{Def embedding} is an isomorphism between $L^{2}(G)$ and $L^{2}(\Gamma)$.

\item
The embedding \eqref{Def embedding} is an isomorphism between $H^{1}(G)$ and $H^{1}(\Gamma)$.

\item
The embedding \eqref{Def embedding} is an isomorphism between $H_{\,0}^{1}(G)$ and $H^{1}_{\,0}(\Gamma)$.
\end{enumerate}
\begin{proof} By definition $\jmath: G\rightarrow \Gamma$ is linear and bijective, therefore the map \eqref{Def embedding} is bijective between spaces of functions. 
\begin{enumerate}[(i)]

\item
Due to the hypothesis $\zeta$ satisfies $C_{1} = \essinf\{\n(s)\cdot\kversor: s\in \Gamma\} > 0$ then, for any $\phi \in L^{2}(\Gamma)$ 
\begin{equation*}
\int_{\Gamma}  (\phi\circ \jmath)^{2} d S = 
\int_{G} (\n(\xthilde)\cdot\kversor)^{-1}  \phi^{2}(\xthilde)\,d\xthilde\leq
\frac{1}{C_{1}}\int_{G}  \phi^{2}(\xthilde)\,d\xthilde .
\end{equation*}
The inequality above gives the continuity of the application $\varphi\mapsto \varphi\circ j$. Due to Banach's inversion theorem the map is an isomorphism.   

\item
By definition $\gradth (\phi \circ \jmath)  = \gradth \phi (\xthilde)$ holds for any $\phi\in H^{1}(G)$.

\item
Is immediate from (ii).
\end{enumerate}
\end{proof}
\end{theorem}
Fix $i\in \{1, \ldots, I\}$, choose $q\in Q_{0}$ supported inside $\Omega(\hi, \Gammai)$ and test  equation \eqref{Pblm lower dimensional variational statement 2}; hence
\begin{multline}\label{Eq tested divergence}
-\int_{\Omega(\hi, \Gammai)} \!\! \utan\cdot\gradth q 
- \int_{\Gammai} \! (\uone\cdot\ni\,)\,\qtwo\, d S
+ \int_{\Gammai+\hi} \!\! (\uone\cdot\ni\,)\,\qtwo\, d S
= \int_{\Gammai\cup \Gammai+\hi} \!\! f_{\scriptscriptstyle \Gammai} \qtwo\,d S .
\end{multline}
We focus on the first term of the left hand side. First $\partial_{z}\qtwo = 0$ implies $\utan\cdot\gradth \qtwo = \utwo\cdot\grad \qtwo$, then
\begin{equation}\label{Eq integration by parts}
- \int_{\Omega(\hi, \Gammai)}\utwo\cdot\grad q 
= \int_{\Omega(\hi, \Gammai)} \div \utwo \,\qtwo  -
\int_{\partial\Omega(\hi, \Gammai)} \qtwo\,\utwo\cdot \outer^{(i)}\; d S .
\end{equation}
The two summands of the left hand side are treated separately. For the first summand the independence from the variable $z$ implies $\div \utwo = \gradth \cdot \utan$, the fact $\utwo\cdot\ni = 0$ in $\Omega(\hi, \Gammai)$ gives $\utan = \MTtang \utang$. Thus
\begin{equation}\label{Eq reduction of divergence term}
\int_{\Omega(\hi, \Gammai)} \!\! \gradth \cdot \utan \,\qtwo\, d\x 
=\hi \int_{\Gi} \!\! \gradth\cdot (\MTtang\utang) \,\qtwo\, d \xthilde
=\hi \int_{\Gammai} \!\!  (\ni\cdot\kversor) \gradth\cdot (\MTtang\utang) \,\qtwo\, d S .
\end{equation}
The boundary term in \eqref{Eq integration by parts} can be written as
\begin{equation*}
-\int_{\Gammai\cup\,\Gammai+\hi} \qtwo\,\utwo\cdot \outer^{(i)}\, d S  
-\int_{\partial\Omega(\hi, \Gammai) - \{\Gammai\cup\,\Gammai+\hi\}} \qtwo\,\utwo\cdot \outer^{(i)}\, d S .
\end{equation*}
The first summand vanishes since $\utwo\cdot \ni = 0$ in $\Omega(\hi, \Gammai)$. The boundary piece described in the second summand is a vertical wall, then $\outer^{(i)}\cdot\kversor = 0$ and it can be identified with the outwards normal vector to the set $\Gi\subseteq\R^{\,2}$. Moreover, due to the independence of the integrand with respect to the variable $z$, the surface integral can be collapsed to a line integral over $\partial\,\Gi$. Combining these observations with \eqref{Eq reduction of divergence term} and \eqref{Eq integration by parts} the equation \eqref{Eq tested divergence} transforms in 
\begin{multline*}
\hi \int_{\Gammai} (\ni\cdot\kversor)\,\gradth\cdot (\MTtang\utang) \,\qtwo\, d S 
- \hi \int_{\partial\Gi} \qtwo\,\u\cdot \outer^{(i)}\, d C\\
- \int_{\Gammai}(\uone\cdot\ni\,)\,\qtwo\, d S
+ \int_{\Gammai+\hi}(\uone\cdot\ni\,)\,\qtwo\, d S
= \int_{\Gammai\cup\,\Gammai+\hi}f_{\scriptscriptstyle \Gamma}\,\qtwo\,dS
\end{multline*}
Where $d C$ is the arc-length measure on $\partial \Gi$. The isomorphisms provided by theorem \eqref{Th manifold isomorphism} imply that the quantifier $\qtwo\vert_{\Gammai} $ can hit any function in the space $ H_{0}^{1}(\Gammai)$.  Therefore, the equation \eqref{Pblm strong limit divergence Omega 2} follows. Finally, using again theorem \eqref{Th manifold isomorphism} the trace of test function $\qtwo\vert_{\Gammai} $ can hit any function in the space $ H_{0}^{1}(\Gammai)$ and combined with equation \eqref{Pblm strong limit divergence Omega 2} give \eqref{Pblm strong limit flux boundary condition}.
%
%
%
%
%
%
\subsection{Strong Convergence of the Solutions}\label{Sec strong convergence}
\begin{theorem}\label{Th strong convergence of solutions}
   Under the hypothesis 
\begin{align}\label{Hyp strong convergence of forcing terms}
& \Vert\,F^{\,\epsilon,\,1}-F^{\,1}\,\Vert_{\,0,\,\Omega_{1}}\rightarrow 0\,, &  
&\Vert\,f^{\,\epsilon}_{\scriptscriptstyle \Gamma} -
f_{\scriptscriptstyle \Gamma}\,\Vert_{\,0,\,\Gamma}\rightarrow 0 \quad \text{and} &
& \Vert\,\g^{\,\epsilon}- \g\,\Vert_{\,0,\,\Omega}\rightarrow 0
\end{align}
The solutions $\ueps, \peps$ satisfy the following strong convergence statements
\begin{align}\label{Hyp strong convergence of solutions}
& \Vert\,\uepsone - \uone\,\Vert_{0, \Omega_{1}}\rightarrow 0 \, , &  
&\Vert\,\epsilon\,\uepstwo - \utwo \Vert_{0 , \Omega_{2}}\rightarrow 0 \, , \\
& \Vert\,\pepsone - \pone\Vert_{1, \Omega_{1}}\rightarrow 0 \, , &
& \Vert\,\pepstwo - \ptwo\Vert_{1, \Omega_{2}}\rightarrow 0 .
\end{align}
\begin{proof}
The proof uses exactly the same arguments presented in \cite{MoralesShow2}, theorem 3.2.
   \end{proof}
\end{theorem}
Finally, assume that $\utwo_{\,\tau}\neq 0$ and consider the quotients:
\begin{equation}\label{Eq blow up relation}
\frac{\left\Vert\,\utangeps\,\,\right\Vert_{\,0,\,\Omega_2}}{\left\Vert\,\uepstwo\cdot\n\,\right\Vert_{\,0,\,\Omega_2}}=
\frac{\left\Vert\,\epsilon\,\utangeps\,\right\Vert_{\,0,\,\Omega_2}}{\left\Vert\,\epsilon\,\uepstwo\cdot\n\,\right\Vert_{\,0,\,\Omega_2}}>
\frac{\left\Vert\,\utwo_{\,\tau}\,\right\Vert_{\,0,\,\Omega_2}-\delta}{\left\Vert\,\epsilon\,\uepstwo\cdot\n\,\right\Vert_{\,0,\,\Omega_2}}>0
\end{equation}
The lower bound holds true for $\epsilon>0$ small enough and adequate $\delta>0$ then we conclude that the magnitudes ratio of the flux tangential component over normal component blows-up to infinity, i.e. the flow in the thin channel is predominantly tangential. Finally if $\utwo_{\tau}=0$, unlike the analysis for flat interfaces presented in \cite{MoralesShow2}, no conclusions can be obtained due to the complexity introduced by the geometry of the fissures.
\begin{figure}[!]
\caption[2]{System of 2-D Manifold Fissures}\label{Fig Lower dimensional Fissured Medium}
\centerline{\resizebox{10cm}{7cm}
{\includegraphics{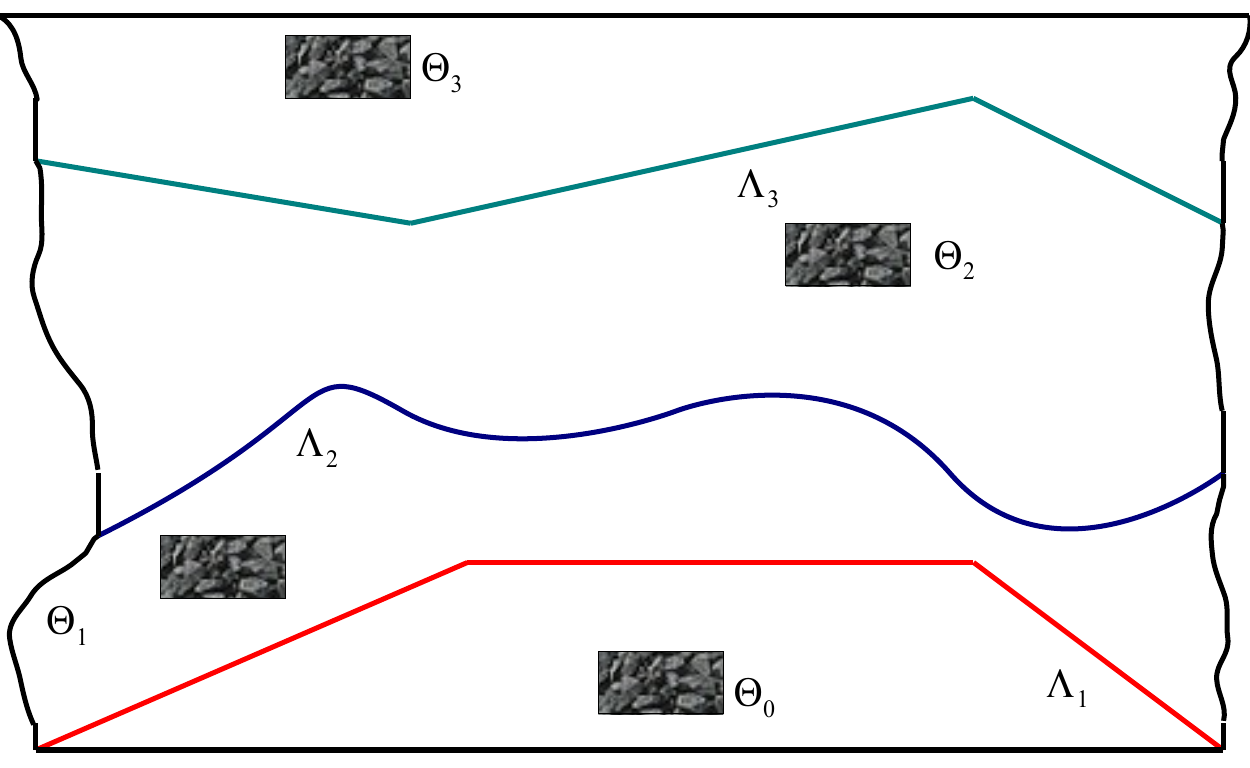} } }
\end{figure}
%
%
\section{A Problem with two dimensional Manifolds}\label{Sec lower dimensional mixed formulation}
In this section, using the independence of the limit functions with respect to $z$ in $\Omega_{2}$ it will be shown that the limiting problem \eqref{Pblm limit problem strong form} can be formulated as a system coupling Darcy flow in three dimensions with tangential flow hosted two dimensional manifolds as depicted in figure \eqref{Fig Lower dimensional Fissured Medium}. 
%
%
%
\subsection{Geometric Setting}\label{Sec limit geometric setting}
\begin{definition}\label{Def lower dimensional medium}
We say a totally fractured medium of two dimensional manifold fissures is a finite collection of
\begin{subequations}\label{Def lower dimensional domain}

Surface functions
\begin{multline}\label{Def fissures surfaces l-d}
\{\lambdai\in C(\overline{\Gi} ):\Gi\subseteq\R^{\,2}\;\text{open bounded simply connected region};\\
\lambdai\;\,\text{piecewise} \;C^{1}
\text{function such that}\;\essinf \ni\cdot\kversor >0\,, \,\iindex \} .
\end{multline}

And rock-matrix regions
\begin{equation}\label{Def rock matrix l-d}
\left\{\Thetaj \subseteq \R^{3}: \Thetaj \neq \emptyset\;
\text{open bounded simply connected region},\; \jindex\right\} .
\end{equation}
%
%
\end{subequations}
Verifying the following properties
\begin{subequations}\label{Def properties of the medium l-d}

Non-overlapping condition and indexed ordered
\begin{equation}\label{Def ordering condition 1 l-d}
\sup \left\{\lambdai(\xthilde): \xthilde\in \Gi\right\} <
\inf \left\{\lambda_{\,i+1}(\xthilde): \xthilde\in G_{i+1}\right\} 
,\;\forall\; 1\leq i\leq I-1
\end{equation}

The interface domain condition 
\begin{equation}\label{Def interface boundary conditions l-d}
\Lambdai = \partial\Thetai \cap \partial\Theta_{i-1}\,,\quad
\forall \,\iindex
\end{equation}

for $\Lambdai \defining \{[\xthilde, \lambdai(\xthilde): \xthilde \in \Gi]\}$. And the connectivity through fissures condition
\begin{equation}\label{Def fissures connectedness condition l-d}
cl(\Theta_{\,\ell})\cap cl(\Theta_{\,k}) = \emptyset\;\,\text{whenever}\;\,
\vert \ell - k\vert > 1 .
\end{equation}
\end{subequations}

For convenience of notation define $\Lambda_{0} \defining \partial \Omega_{0} - \Gamma_{1}$. We denote this fissured system by $\left\{\left(\Lambdai, \Thetai\right):0\leq i\leq I\right\}$. The rock matrix and fissures regions are the sets
\begin{equation}\label{Def fissured medium region l-d}
\begin{split}
\Theta\defining\bigcup_{i = 0} ^{I}\Thetai\,&,\quad
\Lambda\defining \bigcup_{i = 0}^{I} \Lambdai \, ,\\
\Theta_{\scriptscriptstyle FR} & \defining  \Theta\cup\Lambda .
\end{split}
\end{equation}
$\ni$ indicated the upwards normal vector to the surface $\Lambdai$. Finally,  
we introduce the notations $\Lambdaip$ and $\Lambdain$ for the upper and lower faces of the manifold $\Lambdai$.
\end{definition}
%
%
%
\subsection{Spaces of Functions and Isomorphisms}\label{Sec l-d spaces of functions and morphisms}
\begin{definition}\label{Def l-d fissured system spaces}
We define the following spaces for velocity and pressure
\begin{subequations}\label{Eq l-d fissured system spaces}
\begin{multline}\label{Def l-d fissured space of velocities}
\V_{\!\! f} \defining \{\v\in \Ltwo(\Theta_{\scriptscriptstyle FR}):\div \vone\in \Ltwo(\Thetaj),\,\jindex ;\\
\vone\cdot\ni\vert_{\Lambdaip},\, \vone\cdot\ni\vert_{\Lambdain}\,\in L^{2}(\Lambdai)\,,
\vtwo\vert_{\Lambdai}\in \Ltwo(\Lambdai),\,\iindex
\} ,
\end{multline}
\begin{equation}\label{Def l-d fissured space of pressures}
Q_{\! f} \defining \{q\in L^{2}(\Theta_{\scriptscriptstyle FR}):  q\vert_{\Lambdai}\in H^{1}(\Lambdai),\,\iindex\} .
\end{equation}
Endowed with the norms coming from the natural inner products
\begin{multline}\label{Def l-d norm fissured space of velocities}
\Vert\v\Vert_{\scriptscriptstyle \V_{f}} \defining \{\Vert \v \Vert_{\Ltwo(\Theta_{\scriptscriptstyle FR})}^{\,2}
+\Vert \div \v \Vert_{L^{2}(\Theta_{\scriptscriptstyle FR})}^{\,2}\\
+\sum_{ i }\Vert \v\cdot\ni\vert_{\Lambdaip} \Vert_{L^{2}(\Lambdai)}^{\,2}+
\Vert \v\cdot\ni\vert_{\Lambdain} \Vert_{L^{2}(\Lambdai)}^{\,2}+
\Vert \v \Vert_{\Ltwo(\Lambdai)}^{\,2}\}^{1/2} ,
\end{multline}
\begin{equation}\label{Def l-d norm fissured space of pressures}
\Vert q \Vert_{Q_{f}} \defining \{\Vert  q \Vert_{L^{2}(\Theta_{\scriptscriptstyle FR})}^{\,2}
+\sum_{i}\Vert q \Vert_{H^{1}(\Lambdai)}^{\,2}\}^{1/2} .
\end{equation}
\end{subequations}
\end{definition}
\begin{remark}\label{Rem structure of the collapsed space}
Notice that definition \eqref{Def l-d fissured space of velocities} demands only $\v\in \Hdiv(\Thetai)$ \emph{i.e.} the divergence is square integrable only on these subdomains. Therefore, both normal traces $\v\cdot\ni\vert_{\Lambdaip}$ and $\v\cdot\ni\vert_{\Lambdain}$ make sense in $H^{-1/2}(\Gammai)$ but we require the extra condition of been in $L^{2}(\Lambdai)$. We do not demand the global condition $\v\in \Hdiv(\Theta_{\scriptscriptstyle FR})$ because this would imply the continuity of the normal traces across a surface \emph{i.e.} $\uone\cdot\ni\vert_{\Lambdaip} =  \uone\cdot\ni\vert_{\Lambdain}$. Such condition can not model jumps across the fissures as the normal stress balance interface \eqref{Pblm strong normal stress balance} and the limit equation \eqref{Pblm strong limit divergence Omega 2}.
\end{remark}
Next define a change of variable based on piecewise translations
\begin{definition}\label{Def translation change of variable}
Let $\x = (\xthilde, x_{3}) $ and define the map $T:\Omega\rightarrow\R^{\,3}$
\begin{equation}\label{Eq translation change of variable}
T \x\defining \sum_{j \, = \, 0}^{I}\left(\xthilde, x_{3} - \sum_{\ell\, = \,0}^{\,j} h_{\,\ell}\right)\ind_{\Omegaj}(\xthilde, \xthree)
- \sum_{i \, = \, 1}^{I}
\left(\xthilde, \zetai(\xthilde) - 
\sum_{\ell \, = \, 0}^{\,i-1} h_{\,\ell}\right) \ind_{\Omega(\hi, \Gammai)}(\xthilde, \xthree)
\end{equation}
Define $\Thetaj \defining T(\Omegaj)$ and $\lambdai: \Gi\rightarrow\R$ by $\lambdai\defining \zetai(\xthilde) - \sum_{\,\ell = 0}^{\,i-1} h_{\,\ell}$.
\end{definition}
Clearly the system $\{(\Lambdai, \Thetai):\iindex\}$ satisfies the conditions of definition \eqref{Def lower dimensional medium}. With the previous definitions we have the following result
\begin{theorem}\label{Th l-d functions space isomorphisms}
\begin{enumerate}[(i)]
\item
The application $\v \mapsto \v\circ T$ is an isometric isomorphism from $\V_{\! 0}$ to $\V_{\!\!f}$.

\item
The application $q\mapsto q\circ T$ is an isometric isomorphism from $Q_{0}$ to $Q_{\! f}$.
\end{enumerate}
\begin{proof} 
\begin{enumerate}[(i)]
          \item 
          The proof is a direct application of part (i) in theorem \eqref{Th manifold isomorphism}. The only detail that needs further clarification is to observe that
          \begin{equation*}
\begin{split}
\vone\cdot\ni\vert_{\Gammai}\mapsto (\vone\circ T)\cdot\ni\vert_{\Lambdain} \\
\vone\cdot\ni\vert_{\Gammai+\hi}\mapsto (\vone\circ T)\cdot\ni\vert_{\Lambdaip}
\end{split}
\end{equation*}

\item
It is a direct application of parts (i) and (ii) in theorem \eqref{Th manifold isomorphism}.
\end{enumerate}
\end{proof}
\end{theorem}
%
%
\subsection{The Lower Dimensional Mixed Problem}\label{Sec l-d mixed problem}
Due to the previous theorem the problem \eqref{Pblm limit system} is equivalent to the following mixed problem with two dimensional manifolds
\begin{subequations}\label{Pblm reduced dimension variational statement}
\begin{equation*}
 \text { Find } p\in Q_{\! f}\,,\; \u\in \V_{\!\! f} :
\end{equation*}
\begin{multline}\label{Pblm reduced dimension variational statement 1}
\int_{\Theta} a_1 \,\u \cdot \v 
- \int_{\Theta} p \,\div\v  
+  \sum_{i}\hi\int_{\Lambdai} (\ni\cdot\kversor\,)(a_2\, \utang + (\MTtang_{i})'\,
\widetilde{\grad} p +\g_{\,\tau})\cdot \vtang \,d S\\
+\sum_{ i }\alpha\int_{\Lambdai}\left[(\uone\cdot\ni\vert_{\Lambdaip})(\vone\cdot\ni\vert_{\Lambdaip}) +(\uone\cdot\ni\vert_{\Lambdain})(\vone\cdot\ni\vert_{\Lambdain})\right] d S\\
-\sum_{ i }\int_{\Lambdai}\ptwo \left[(\vone\cdot\ni\vert_{\Lambdaip}) - (\vone\cdot\ni\vert_{\Lambdain})\right] d S
= - \int_{\Theta}\g \cdot \v - \int_{\Lambda}\g_{\tau} \cdot \vtang \, dS
\end{multline}
\begin{multline}\label{Pblm reduced dimension variational statement 2}
\int_{\Theta}\div\u\, q\, 
- \sum_{i} \hi \int_{\Lambdai} (\ni\cdot\kversor\,)\,
\MTtang_{i}\utang \cdot \gradth \qtwo \,d S\\
+ \sum_{ i }\int_{\Lambdai}\left[(\uone\cdot\n\vert_{\Lambdaip}) - (\uone\cdot\n\vert_{\Lambdain})\right]\,\qtwo \, d S
= \int_{\Theta}F^{\, 1}\, q 
+ \int_{\Lambda} f_{\scriptscriptstyle \Gamma}\, \qtwo\,d S
\end{multline}
\begin{equation*}
\text{for all}\;q\in Q_{\!f}\,,\; \v\in \V_{\!\! f} .
\end{equation*}
\end{subequations}
Finally the equivalence of problems \eqref{Pblm limit system} and \eqref{Pblm lower dimensional variational statement}  gives the well-posedness of the system above.
%
%
%
%
%
%
\section{Final Discussion and Future Work}\label{Sec limitations and future work}
\begin{figure}[b]
\caption[3]{Translation Generated Fissures}\label{Fig Fissure with high gradients}
\hspace{8pt}
\centerline{\resizebox{12cm}{5cm}
{\includegraphics{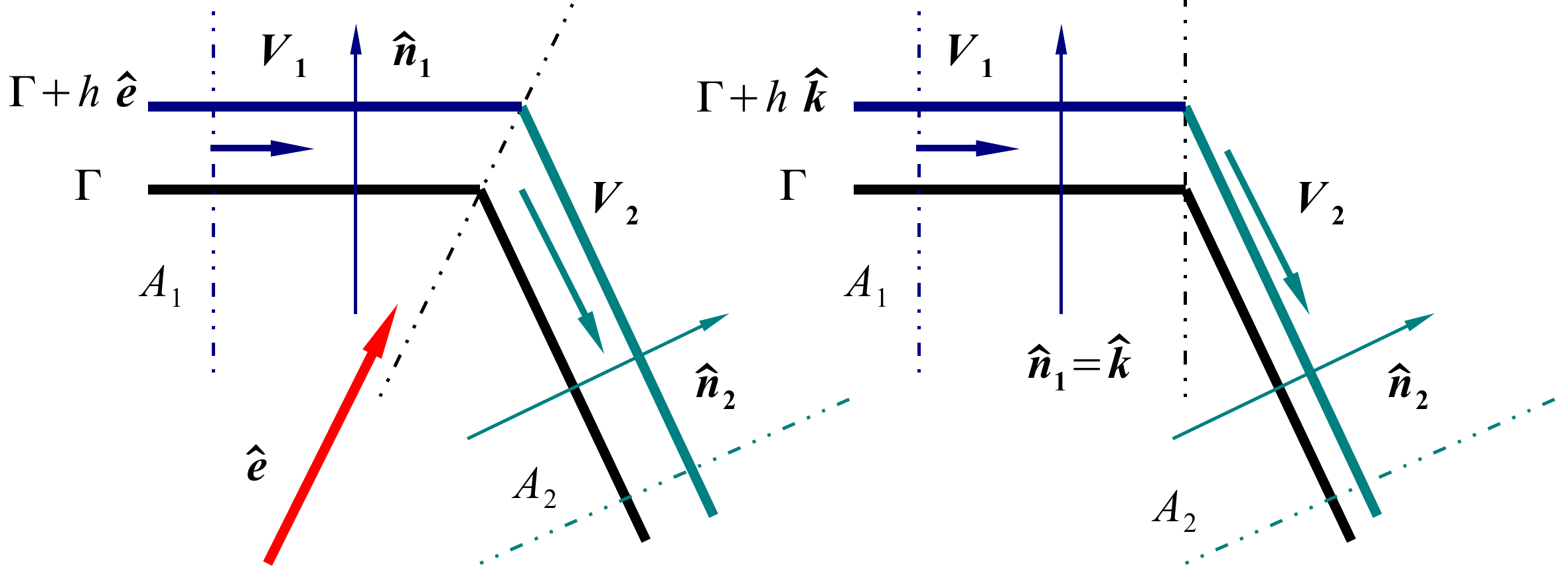} } }
\end{figure}
\begin{enumerate}[(i)]

\item
The formulation presented in this work can manage large amounts of information in a remarkably efficient way. One of the main reasons is the notation introduced by Showalter in \cite{MoralesShow2} for the description of function spaces.  

\item
The results can be generalized immediately to the $\R^{\! N}$-setting using the same arguments presented here. The structure of the problems is analogous.

\item
The approach based on analytic semigroups theory presented in section \cite{MoralesShow2} can be directly applied here to model the time dependent problem for totally fissured systems with singularities. 

\item
Although the mathematical analysis is solid, the approach used throughout the paper stops been suitable for surfaces with high gradients such as the one depicted in the right hand side of figure \eqref{Fig Fissure with high gradients} where $\n_{2}\cdot\kversor \ll \n_{1}\cdot\kversor$. In this case the translation in the direction $\kversor$ generates a fissure whose cross section areas can be very different from one piece to another i.e. $A_{2}\ll A_{1}$. Such a fissure is not realistic. On the other hand consider a fissure such as the one depicted in the left hand side of figure \eqref{Fig Fissure with high gradients}. Here the translation is made in the bisector vector direction 
\begin{equation*}
\eversor = \frac{1}{\vert \frac{\n_{1}+\n_{2}}{2} \vert}\,\frac{\n_{1}+\n_{2}}{2}  
\end{equation*}
This process generates a more realistic fissure. Additionally, demanding the fissures to be defined by the parallel translation of a surface in a fixed direction although is a step forward with respect to previous achievements, is still restrictive for modeling the phenomenon in natural geological formations. Setting the problem in the mixed variational formulation used here can be easily extended to systems with fissures described by a very general type of geometry. However, the difficulty of the asymptotic analysis increases substantially. 

Such question will be addressed in future work by the introduction of correction factors obtained comparing the flow energy dissipation in a real fissure and an artificial one e.g. replacing the presence of the fissure in the left hand side of \eqref{Fig Fissure with high gradients} with the one on the right side affected by a correction factor. In the same way, fissures defined by walls which are not rigid translations of the other will be compared to a fissure generated by vertical translation of its ``average surface'' and having the same ``average width''.
\end{enumerate}
%
%
%
%
\section{Acknowledgements}\label{Sec Acknowledgements}
The author thanks to Universidad Nacional de Colombia, Sede Medell\'in for partially supporting this work under the projects HERMES 17194 and HERMES 14917 as well as the Department of Energy, Office of Science, USA for partially supporting this work under grant 98089. Finally, the author wishes to thank Professor Ralph Showalter from the Mathematics Department at Oregon State University for his helpful insight, observations and suggestions.  
\nocite{Showalter77}
%
%
%
%


\begin{thebibliography}{10}

\bibitem{Allaire2009}
Gr\'egorie Allaire, Marc Briane, Robert Brizzi, and Yves Capdeboscq.
\newblock Two asymptotic models for arrays of underground waste containers.
\newblock {\em Applied Analysis}, 88 (no. 10-11):1445--1467, 2009.

\bibitem{ArbBrunson2007}
Todd Arbogast and Dana Brunson.
\newblock A computational method for approximating a {D}arcy-{S}tokes system
  governing a vuggy porous medium.
\newblock {\em Computational Geosciences}, 11, No 3:207--218, 2007.

\bibitem{ArbLehr2006}
Todd Arbogast and Heather Lehr.
\newblock Homogenization of a darcy-stokes system modeling vuggy porous media.
\newblock {\em Computational Geosciences}, 10, No 3:291--302, 2006.

\bibitem{Gunzburger2009}
Nan Chen, Max Gunzburger, and Xiaoming Wang.
\newblock Asymptotic analysis of the differences between the {S}tokes-{D}arcy
  system with different interface conditions and the {S}tokes-{B}rinkman
  system.
\newblock {\em Journal of Mathematical Analysis and Applications}, 368
  (2):658--676, 2009.

\bibitem{Gatica2009}
Gabriel~N. Gatica, Salim Meddahi, and Ricardo Oyarz\'ua.
\newblock A conforming mixed finite-element method for the coupling of fluid
  flow with porous media flow.
\newblock {\em IMA Journal of Numerical Analysis}, 29, 1:86--108, 2009.

\bibitem{GiraultRaviartFEM}
V.~Girault and P.-A. Raviart.
\newblock {\em Finite element approximation of the {N}avier-{S}tokes
  equations}, volume 749 of {\em Lecture Notes in Mathematics}.
\newblock Springer-Verlag, Berlin, 1979.

\bibitem{Hornung}
Ulrich Hornung, editor.
\newblock {\em Homogenization and Porous Media, Ulrich Hornung editor},
  volume~6 of {\em Interdisciplinary Applied Mathematics}.
\newblock Springer-Verlag, New York, 1997.

\bibitem{ZhaoQin}
ZhaoQin Huang, Jun Yao, YaJun Li, ChenChen Wang, and XinRui L\"{u}.
\newblock Permeability analysis of fractured vuggy porous media based on
  homogenization theory.
\newblock {\em Science China Technological Sciences}, 53 (3):839--847, 2010.

\bibitem{Yotov}
W.~J. Layton, F.~Scheiweck, and I.~Yotov.
\newblock Coupling fluid flow with porous media flow.
\newblock {\em SIAM Journal of Numerical Analysis}, 40 (6):2195--2218, 2003.

\bibitem{Levy83}
Th{\'e}r{\`e}se L{\'e}vy.
\newblock Fluid flow through an array of fixed particles.
\newblock {\em International Journal of Engineering Science}, 21:11--23, 1983.

\bibitem{Loredana1}
J.~San Mart{\'i}n, J.-F. Scheid, and L.~Smaranda.
\newblock A modified lagrange-galerkin method for a fluid-rigid system with
  discontinuous density.
\newblock {\em Numerische Mathematik}, 122 (2):341--382, 2012.

\bibitem{JaffRob05}
Vincent Martin, J{\'e}r{\^o}me Jaffr{\'e}, and Jean~E. Roberts.
\newblock Modeling fractures and barriers as interfaces for flow in porous
  media.
\newblock {\em SIAM J. Sci. Comput.}, 26(5):1667--1691, 2005.

\bibitem{Mikelic89}
A.~Mikeli{\'c}.
\newblock A convergence theorem for homogenization of two-phase miscible flow
  through fractured reservoirs with uniform fracture distribution.
\newblock {\em Applicable Anal.}, 33:203--214, 1089.

\bibitem{MoralesShow1}
Fernando Morales and Ralph Showalter.
\newblock The narrow fracture approximation by channeled flow.
\newblock {\em Journal of Mathematical Analysis and Applications},
  365:320--331, 2010.

\bibitem{MoralesShow2}
Fernando Morales and Ralph Showalter.
\newblock Interface approximation of {D}arcy flow in a narrow channel.
\newblock {\em Mathematical Methods in the Applied Sciences}, 35:182--195,
  2012.

\bibitem{SP80}
Enrique S\'anchez-Palencia.
\newblock {\em Nonhomogeneous media and vibration theory}, volume 127 of {\em
  Lecture Notes in Physics}.
\newblock Springer-Verlag, Berlin, 1980.

\bibitem{Showalter77}
R.~E. Showalter.
\newblock {\em Hilbert space methods for partial differential equations},
  volume~1 of {\em Monographs and Studies in Mathematics}.
\newblock Pitman, London-San Francisco, CA-Melbourne, 1977.

\bibitem{Showalter97}
R.~E. Showalter.
\newblock {\em Monotone operators in {B}anach space and nonlinear partial
  differential equations}, volume~49 of {\em Mathematical Surveys and
  Monographs}.
\newblock American Mathematical Society, Providence, RI, 1997.

\bibitem{Show97}
R.E. Showalter.
\newblock {\em Microstructure Models of Porous Media. In Ulrich Hornung editor
  Homogenization and Porous Media}, volume~6 of {\em Interdisciplinary Applied
  Mathematics}.
\newblock Springer-Verlag, New York, 1997.

\end{thebibliography}
\end{document}